\documentclass[review]{elsarticle}

\usepackage{lineno,hyperref}

\usepackage{amsmath}
\usepackage{amssymb}
\usepackage{amsthm}
\usepackage{mathrsfs}
\usepackage{graphicx}
\usepackage{epstopdf}
\usepackage{bm}
\usepackage{color}
\usepackage{catoptions}
\usepackage{hyperref}
\usepackage{cleveref}[2011/03/22]

\Crefformat{figure}{#2Fig.~#1#3}
\Crefformat{equation}{#2Eq.~(#1#3)}

\makeatletter
\def\setupcrefformat{%
    \docommalist{algorithm,appendix,chapter,corollary,definition,
    enumi,Equation,Example,Figure,footnote,lemma,line,note,part,
    proposition,remark,result,section,Table,theorem} 
    {%
    \crefformat{##1}{####2##1~####1####3}%
     We don't seem to need the following:
     \Crefformat{##1}{####2\@nameuse{Cref@##1@name}~####1####3}%
  }%
}
\makeatother

\makeatletter
\renewcommand{\fnum@figure}{Fig. \thefigure}
\makeatother

\usepackage[numbers]{natbib}

\modulolinenumbers[10]
\theoremstyle{definition}
\newtheorem{exmp}{Example}[subsection]
\numberwithin{equation}{section}

\usepackage{geometry}
\makeatletter
\newcommand{\mathleft}{\@fleqntrue\@mathmargin0pt}
\newcommand{\mathcenter}{\@fleqnfalse}
\makeatother

\journal{Journal of \LaTeX\ Templates}
\makeatletter
\def\ps@pprintTitle{%
    \let\@oddhead\@empty
    \let\@evenhead\@empty
    \let\@oddfoot\@empty
    \let\@evenfoot\@oddfoot
}
\makeatother

\begin{document}

\begin{frontmatter}

\title{High order residual distribution for steady state problems for hyperbolic conservation laws}

\author{Jianfang Lin\fnref{J.Lin}}
\fntext[J.Lin]{School of Mathematical Sciences, Xiamen University, Xiamen, Fujian 361005, PR China. Email: jflin@stu.xmu.edu.cn}

\author{R\'emi Abgrall\fnref{R.Abgrall}}
\fntext[R.Abgrall]{Institute of Mathematics, University of Zurich, Zurich 8057, Switzerland. Email: remi.abgrall@math.uzh.ch}

\author{Jianxian Qiu\fnref{J.Qiu}}
\fntext[J.Qiu]{School of Mathematical Sciences and Fujian Provincial Key Laboratory of Mathematical Modeling and High-Performance Scientific Computing, Xiamen University, Xiamen, Fujian 361005, PR China. Email: jxqiu@xmu.edu.cn}

%




\begin{abstract}

In this paper, we propose a high order residual distribution conservative finite difference scheme for solving steady state conservation laws. A new type of WENO (weighted essentially non-oscillatory) termed as WENO-ZQ integration is used to compute the numerical fluxes and source term based on the point values of the solution, and the principles of residual distribution schemes are adapted to obtain steady state solutions. Extensive numerical examples in both scalar and system test problems in one and two dimensions demonstrate the efficiency, high order accuracy and the capability of resolving shocks of the proposed methods.

\end{abstract}


\end{frontmatter}


\section{Introduction}

We consider hyperbolic conservation laws with source terms
\begin{equation}\label{Eq:hcl}
u_{t} + \nabla \cdot f (u) = s(u, x),
\end{equation}
in which the Jacobian matrix $f'(u)$ is diagonalizable with all the eigenvalues being real for any u. In recent decades many high order methods, such as finite difference methods, finite volume methods and discontinuous Galerkin (DG) methods, have been investigated to solve for hyperbolic conservation laws. Within these schemes, the essentially non-oscillatory (ENO) and weighted ENO (WENO) reconstructions \cite{Harten.Osher_JCP1987, Shu.Osher_JCP1988, Shu.Osher_JCP1989, Liu.Osher.Chan_JCP1994, Jiang.Shu_JCP1996, Shu_Note_1998} are very successful in capturing shocks in a sharp, non-oscillatory fashion while maintaining high order accuracy in smooth regions. Later on, the various types of ENO and WENO schemes are quite successful in numerical simulations for steady state and unsteady problems harboring strong discontinuities and sophisticated smooth structures. Recently, a new type of WENO termed as WENO-ZQ schemes \citep{Zhu.Qiu_JCP2016, Zhu.Qiu_JSC2017} was proposed, which has the advantages of simplicity, high order accuracy and easy implementation in the computation.

In this paper, we are interested in computing the steady solution of \eqref{Eq:hcl} and developing high order conservative schemes which are of finite difference type (the numerical approximations are the point values of the solution) and have a comparable computational cost as regular finite difference schemes, the meshes are allowed to be arbitrary Cartesian or curvilinear without any smoothness assumption. Many current schemes use ideas for high resolution schemes developed in the 1970s and 1980s by van Leer, Roe, Osher, Harten, Yee, Sweby and many others \citep{Harten.Lax.vanLeer_SIAM1983, Roe_Springer1987, Engquist.Osher_MC1980, Osher.Solomon_MC1982, Osher.Chakravarthy_SIAM1984, Osher.Chakravarthy_Springer1986, Harten.Lax_SIAM1981, Yee_CMA1986, Sweby_SIAM1984}. However, the quality of the solution is still questionable: some apparently simple problems, such as computing the lift and drag of an airfoil, still pose difficulties. One reason is that the so-called high resolution schemes suffer a much too great entropy production. In fact, they have been on one dimensional scalar problems, then extended to multi-Dimension systems, but their construction relies on ``1D ideas". Another difficult problem is the sensitivity to the mesh. It is still difficult to construct a 3D mesh of consequently, the quality of the solution itself may be questionable in many cases. Hence, it is natural to construct methods that have as little sensitivity as possible to the regularity of mesh.

For these reasons, decades years many researchers have tried to incorporate ideas contained in the 1D high-resolution schemes (upwind) into a finite-element-like framework. Some of the major contributions \citep{Struijs.Deconinck.Roe_CFD1991, Roe.Sidilkover_JNA1992, Deconinck.Struijs.Bourgeois.Roe_CFD1993} have been made by P. L. Roe, H. Deconinck, D. Sidilkover and their coauthors. These residual distribution (RD) or fluctuation splitting schemes, were first developed for a scalar transport equation, then formally extended to systems. These schemes share many common features with the streamline upwind Petrov Galerkin (SUPG) schemes of Hughes \citep{Hughes_CMAME1986} or the streamline diffusion methods of Johnson \citep{Johnson_MC1986}, except for up-winding. A brief view of a RD scheme for  \eqref{Eq:hcl} is given as follows: an approximate solution of \eqref{Eq:hcl} is sought on a general triangular or quadrilateral mesh $\mathscr{T}_{h}$. The nodes of $\mathscr{T}_{h}$ are denoted by $\left\{M_{i}\right\}$ and $T$ is a generic element. On each element $T$, we define a total residual $\Phi^{T}$ and also define $\Phi^{T}_{i}$ as the amount of $\Phi^{T}$ associated with vertex $M_{i}$, such that a conservation property is satisfied
\begin{equation}
\Phi^{T} = \int_{T}\! (\nabla\cdot f^{h}(u_{h}) - s^{h}(u_{h}, x))\, \mathrm{d}x ,~~\sum\limits_{i, M_{i}\in T}\Phi^{T}_{i} = \Phi^{T}.
\end{equation}
Then the residual distribution scheme is given as
\begin{equation}
\left|C_{i}\right|\frac{u^{n+1}_{i}-u^{n}_{i}}{\Delta t_{n}} + \sum\limits_{T, M_{i}\in T}\Phi^{T}_{i} = 0,
\end{equation}
where $\left|C_{i}\right|$ is the area of the dual element associated with $M_{i}$.

Recently, RD schemes have received considerable attention and they are demonstrated to be robust in many  numerical tests. A Lax-Wendroff type theorem has been provided to ensure convergence to the weak solution \citep{Abgrall.Mer.Nkonga_2002}, and the stability is established following maximum principles, see, e.g.\citep{Abgrall_JCP2001, Abgrall.Mezine_JCP2004}. The accuracy at steady state is ensured if the scheme satisfies the residual property, which is related to the accuracy of approximating the residuals, see \citep{Abgrall_JCP2001}. The works mentioned above are mostly for schemes of at most second order accuracy, which follows the systematic construction from a first order monotone and upwind RD scheme to a second one. Later on, RD schemes were generalized to high order schemes by Abgrall and Roe \citep{Abgrall.Roe_JCP2003} on general triangular meshes. Based on the same distribution principles, Chou and Shu \citep{Chou.Shu_JCP2006} developed a finite difference method  based on RD scheme which works on curvilinear meshes, and their scheme achieves high order accuracy and low computational cost as in finite difference methods, but it need add an additional dissipation residual around the shocks, only in two dimensional cases. Motivated by their work, we are interested in developing a finite difference method based on RD scheme and use a SUPG-like distribution properties. Because of the same width of the stencil as the classical WENO reconstruction when using WENO-ZQ integration reconstruction, a Lax-Wendroff theorem for convergence towards weak solutions is the same proof as shown in \citep{Chou.Shu_JCP2006}.

This paper is organized as follows: in the Sections 2 and 3, we describe the residual evaluation and the residual distribution procedures for one and two-dimensional problems, respectively. In the Section 4, the numerical simulation results for one and two-dimensional scalar and system steady state problems are shown to demonstrate the good behaviors of our scheme. Concluding remarks are given in the Section 5.

\section{High order RD finite difference WENO-ZQ schemes in one dimension}

In this section, we design a residual distribution high order WENO-ZQ finite difference scheme for one-dimensional steady state problems. In the first subsection, we define the total residual within each cell from the integral form, and then describe the distribution of the total residual within each cell, complying with the principles of SUPG-like and the residual property. In the second subsection, we extend the scheme naturally to the one-dimensional systems, based on a local characteristic field decomposition, and using the principles as in the scalar case to distribute the total residual within each cell in the characteristic fields.

\subsection{One-dimensional scalar problems}\label{Subsec:scalar1d}

We consider the one-dimensional scalar steady state problem
\begin{equation}\label{Eq:ssp1d}
f(u)_{x} = s(u, x).
\end{equation}
We define the grid to be $\left\{x_{i}\right\}_{i=0,\cdots, N}$, the grid function $\left\{u_{i}\right\}_{i=0, \cdots, N}$, the interval $I_{i+\frac{1}{2}}=\left[x_{i}, x_{i+1}\right]$, the step x-direction $\Delta x_{i+\frac{1}{2}}$, the control volume centered at $x_{i}$ to be $C_{i}$ (from the mid-point of the interval $I_{i-\frac{1}{2}}$ to the mid-point of the interval $I_{i+\frac{1}{2}}$), and the length of $C_{i}$ is denoted by $\left|C_{i}\right|$.

The total residual in the interval $I_{i+\frac{1}{2}}$ is defined by
\begin{equation}\label{Def:residual}
\Phi_{i+\frac{1}{2}}=\int^{x_{i+1}}_{x_{i}}\!(f(u)_{x} - s(u, x))\,\mathrm{d}x = f(u_{i+1}) - f(u_{i}) - \int^{x_{i+1}}_{x_{i}}\! s(u, x)\, \mathrm{d}x.
\end{equation}
If we can reach the zero residual limit, i.e. if $\Phi_{i+\frac{1}{2}}=0$ for all i, the accuracy of the scheme is determined by the accuracy of the approximation to $\int^{x_{i+1}}_{x_{i}}\! s(u,x)\, \mathrm{d}x$. In our scheme, we use a fourth order WENO-ZQ integration to approximate the integral $\int^{x_{i+1}}_{x_{i}}\!s(u, x)\,\mathrm{d}x$ (leading to a fifth order WENO-ZQ approximation to the integral within each cell and hence a fourth order approximation to the integral over the whole computational domain), which is described as follows:

Step 1. Choose the following big stencil: $S_{1} = \left\{x_{i-1}, x_{i}, x_{i+1}, x_{i+2}\right\}$, there is an unique polynomial $p_{1}(x)$ of degree 4 which interpolates $s(u, x)$ at nodes in $S_{1}$ and satisfying:
\begin{equation}
p_{1}(x_{j}) = s(u_{j}, x_{j}), ~~j = i-1, i, i+1, i+2.
\end{equation}
Choose another smaller stencil: $S_{2} =\left\{x_{i}, x_{i+1}\right\}$, there is an unique linear polynomial $p_{2}(x)$ which interpolates $s(u, x)$ at nodes in $S_{2}$ and satisfying:
\begin{equation}
p_{2}(x_{j}) = s(u_{j}, x_{j}), ~~j = i, i+1.
\end{equation}
Then, we integrate $p_{1}(x)$ and $p_{2}(x)$ in the interval $I_{i+\frac{1}{2}}$, denoted by $q_{1}$ and $q_{2}$, respectively.

Step 2. The main selection principle of the linear weight is solely based on the consideration of a balance between the accuracy and the ability to achieve essentially non-oscillatory shock transitions. Here,  we rewrite $q_{1}$ as $q_{1}=\gamma_{1}(\frac{1}{\gamma_{1}}q_{1} - \frac{\gamma_{2}}{\gamma_{1}}q_{2})+ \gamma_{2}q_{2}$. In all of our numerical tests, following the practice in \citep{Dumbser.Kaser_JCP2007, Zhong.Shu_JCP2013}, we take the positive linear weights as $\gamma_{1}=0.99$ and $\gamma_{2}=0.01$. The linear weights can be chosen to be any set of positive numbers on condition that the summation is 1 and would not pollute the new scheme's optimal accuracy.

Step 3. Compute the smoothness indicators $\beta_{n}, n = 1, 2$, which measure how smooth the functions $p_{n}(x), n =1, 2$, are in the target cell $I_{i+\frac{1}{2}}$. The smaller these smoothness indicators, the smoother the functions are in $I_{i+\frac{1}{2}}$. We use the same recipe for the smoothness indicators as in \citep{Balsara.Rumpf_Dumbser.Munz_JCP2009, Jiang.Shu_JCP1996, Shu_SIAMRev2009};
\begin{equation}
\beta_{n} = \sum\limits^{r}_{m=1}\int\limits_{I_{i+\frac{1}{2}}}(\Delta x_{i+\frac{1}{2}})^{2m-1}\!(\frac{d^{m}p_{n}}{dx^{m}})^{2}\, \mathrm{d}x, ~~n=1, 2,
\end{equation}
where $r$ is the degree of the corresponding polynomial.

Step 4. Calculate the non-linear weights based on the linear weights and the smoothness indicators. For instance, as shown in \citep{Borges.Carmona.Costa.Don_JCP2008, Castro.Costa.Don_JCP2011}, we use new $\tau_{0}$ which is simply defined as the square of the absolute difference between $\beta_{1}$ and $\beta_{2}$, and is different to the formula specified in \citep{Borges.Carmona.Costa.Don_JCP2008, Castro.Costa.Don_JCP2011}. It is defined as follows:
\begin{equation}
\tau_{0} = \left|\beta_{1}-\beta_{2}\right|^{2}.
\end{equation}
Then, we define
\begin{equation}
\omega_{n} = \frac{\bar{\omega}_{n}}{\sum^{2}_{l=1}\bar{\omega}_{l}}, ~~ \bar{\omega}_{n} = \gamma_{n}(1+\frac{\tau_{0
}}{\varepsilon + \beta_{n}}), ~~n = 1, 2,
\end{equation}
which satisfy the order accuracy $\omega_{n} = \gamma_{n} + \mathcal{O}(\Delta x_{i +\frac{1}{2}}^{4})$, where $\varepsilon$ is a small positive number to prevent the denominator from becoming zero. And we take $\varepsilon = 10^{-6}$  in our computation.

Step 5. The new final reconstruction of the integral $\int^{x_{i+1}}_{x_{i}}\!s(u, x)\, \mathrm{d}x$ in $I_{i+\frac{1}{2}}$, as shown in \citep{Zhu.Qiu_JCP2016}, is given by
\begin{equation}
\int^{x_{i+1}}_{x_{i}}\!s(u, x)\, \mathrm{d}x = \omega_{1}(\frac{1}{\gamma_{1}}q_{1} - \frac{\gamma_{2}}{\gamma_{1}}q_{2}) + \omega_{2}q_{2} + \mathcal{O}(\Delta x_{i+\frac{1}{2}}^{5}).
\end{equation}

Let us mention that near boundary, one-sided biased rather than central stencils could be used in WENO-ZQ procedure.

Next, we start to distribute the total residuals. In the interval $\left[x_{i}, x_{i+1}\right]$, the total residual is $\Phi_{i+\frac{1}{2}}$, and it is to be distributed to the nodes $x_{i}$ and $x_{i+1}$. For simplicity and with no ambiguity, we drop the subscript $i+\frac{1}{2}$ off for the total residual $\Phi_{i+\frac{1}{2}}$. Here we denote the residuals distributed to the points $x_{i}$ and $x_{i+1}$ as $\Phi^{-}$ and $\Phi^{+}$, respectively. To have a SUPG-like scheme, one way to distribute the total residual $\Phi$ is the following:

Step 1. First order Lax-Friedrichs linear distribution is given by
\begin{equation}\label{Def:LaxF1dscl}
\Phi^{\text{LxF}, -} = \frac{1}{2}\Phi + \alpha (u_{i} - \bar{u}),~~\Phi^{\text{LxF}, +} = \frac{1}{2}\Phi + \alpha (u_{i+1} - \bar{u}),
\end{equation}
where  $\bar{u}$ is an average state in the cell taken to be $\frac{1}{2}(u_{i} + u_{i+1})$, and $\alpha$ is determined by
\begin{equation}
\alpha = \Delta x_{i + \frac{1}{2}}\cdot\max\limits_{j\in I_{i+\frac{1}{2}}}\left\{\left|f'(u_{j})\right|\right\}.
\end{equation}
The Struijs' ``limiter'' is defined in the following:

\begin{eqnarray}
\beta^{-} =  \frac{\max(\Phi^{\text{LxF}, -}/ \Phi, 0)}{\sum\limits_{*\in \{-, +\}}\max (\Phi^{\text{LxF}, *} / \Phi, 0)},\label{Def:Struijs1d-}\\
\beta^{+} =  \frac{\max(\Phi^{\text{LxF}, +}/ \Phi, 0)}{\sum\limits_{*\in \{-, +\}}\max (\Phi^{\text{LxF}, *} / \Phi, 0)}\label{Def:Struijs1d+}.
\end{eqnarray}

Step 2. The streamline dissipation term is defined by
\begin{equation}\label{Def:Streamlinediss}
\int\limits_{I_{i+\frac{1}{2}}}\!(\nabla_{u}f(u)\cdot\nabla\varphi_{j})\tau(\nabla_{u}f(u)\cdot\nabla u -s(u, x))\, \mathrm{d}x, ~~j =i, i+1,
\end{equation}
where $\varphi_{j}$ is the basis function associated to the node j and $\tau > 0$ in the interval $I_{i+\frac{1}{2}}$. We take $\varphi_{i} = -\frac{x-x_{i+1}}{x_{i+1}-x_{i}}$ and $\varphi_{i+1} = \frac{x-x_{i}}{x_{i+1}-x_{i}}$, and $\tau^{-1}$ is defined by
\[
\tau^{-1} = \sum\limits_{j\in \{i, i+1\}}\left|f'(\bar{u})\varphi'_{j}(x_{j})\right|.
\]

As for one dimensional scalar case, we have
\begin{equation}\label{Def:dissipation1d}
\Phi^{-}_{\text{diss}} = -\frac{1}{2}\frac{f'(\bar{u})}{\left| f'(\bar{u})\right|}\Phi,~~ \Phi^{+}_{\text{diss}} = \frac{1}{2}\frac{f'(\bar{u})}{\left| f'(\bar{u})\right|}\Phi.
\end{equation}
In order to prevent the absolute $\left|f'(\bar{u})\right|$ from becoming zero, we take the Roe's correction
\begin{equation}\label{Def:Roecorr}
\left|f'(\bar{u})\right| = \begin{cases}
\hfil \left|f'(\bar{u})\right| & {\texttt{if}}~ \left|f'(\bar{u})\right| > \epsilon,\\
\hfil \frac{f'(\bar{u})^{2}+\epsilon^{2}}{2\epsilon} & \texttt{else},
\end{cases}
\end{equation}
where $\epsilon$ is taken to be $10^{-2}$ in the computation.
 
Hence, we get the way to distribute the total residual within each cell as follows:
\begin{equation}\label{Def:rd1d_scl}
\Phi^{-} =\beta^{-}\Phi + \Phi^{-}_{\text{diss}}, ~~ \Phi^{+} =\beta^{+}\Phi + \Phi^{+}_{\text{diss}}.
\end{equation}

Finally, the point value $u_{i}$ is updated through sending the distributed residuals to the point $x_{i}$, as in a pseudo time-marching scheme, which can be written as a semi-discrete system
\begin{equation}\label{Eq:RDscheme}
\frac{\mathrm{d}u_{i}}{\mathrm{d}t} + \frac{1}{\left|C_{i}\right|}\left(\Phi^{+}_{i-\frac{1}{2}}+\Phi^{-}_{i+\frac{1}{2}}\right)=0.
\end{equation}
In our numerical experiments, we use a third order TVD Runge-Kutta scheme \citep{Shu.Osher_JCP1989} for the (pseudo) time discretization. Since the accuracy in time is irrelevant here, any stable time marching can be used.

We now summarize the procedure of the high order RD finite difference WENO-ZQ scheme for one-dimensional scalar problems:
\vspace{0.2cm}
\\
1. Compute the total residual defined in \Cref{Def:residual} using WENO-ZQ integration with a proper accuracy for the source term.
\\
2. Distribute the total residual within each cell according to the SUPG-like principle, which is defined in \eqref{Def:LaxF1dscl}, \eqref{Def:Struijs1d-}, \eqref{Def:Struijs1d+}, \eqref{Def:dissipation1d}, \eqref{Def:rd1d_scl}.
\\
3. Update the point values through sending the residuals and forward in pseudo time by a TVD Runge-Kutta time discretization until the steady state is reached.

\subsection{One-dimensional systems}

Consider a one-dimensional steady state system \eqref{Eq:ssp1d} where {\bf u}, {\bf f}({\bf u}) and {\bf s}({\bf u}, x) are vector-valued functions in $\mathbb{R}^{m}$. For hyperbolic systems, we assume that the Jacobian ${\bf f}'({\bf u})$ can be written as $R\Lambda L$, where $\Lambda$ is a diagonal matrix with real eigenvalues on the diagonal, and L and R are matrices of left and right eigenvectors of ${\bf f}'({\bf u})$, respectively.

The grid, grid function, step x-direction and control volume
 are denoted as in the \Cref{Subsec:scalar1d}. The total residual ${\bf \Phi}_{i+\frac{1}{2}} $ in the interval $\left[x_{i}, x_{i+1}\right]$ is again defined by \eqref{Def:residual}. As before, the accuracy of the scheme is determined by the accuracy of the approximation to $\int^{x_{i+1}}_{x_{i}}\!{\bf s}({\bf u}, x)\, \mathrm{d}x$, which is again obtained by a fourth order WENO-ZQ integration.

In order to distribute the total residual ${\bf \Phi}_{i+\frac{1}{2}}$, we need use a local characteristic decomposition when we define the Struijs' ``limiter" in the interval $\left[x_{i}, x_{i+1}\right]$. First, we compute an average state $\bar{{\bf u}}$ between ${\bf u}_{i+1}$ and ${\bf u}_{i}$, using either the simple arithmetic mean or Roe's average \citep{Roe_JCP1981}, and $\bar{L}$ and $\bar{R}$ are the corresponding left and right eigenvectors L and R evaluated at the average state $\bar{{\bf u}}$, and $\bar{\lambda}_{k}$ is the corresponding kth eigenvalue. In the following, for simplicity and with no ambiguity, we drop the subscript $i+\frac{1}{2}$ off for the total residual ${\bf \Phi}_{i+\frac{1}{2}}$. The first order Lax-Friedrichs linear distribution is again defined by \eqref{Def:LaxF1dscl}, then we project ${\bf \Phi}^{\text{LxF}, -}$ and ${\bf \Phi}^{\text{LxF}, +}$ to the characteristic fields, namely, ${\bf \Psi}^{\text{LxF}, -}=\bar{L}{\bf \Phi}^{\text{LxF}, -}$ and ${\bf \Psi}^{\text{LxF}, +}=\bar{L}{\bf \Phi}^{\text{LxF}, +}$, respectively, with ${\bf \Psi} = {\bf \Psi}^{\text{LxF}, -} + {\bf \Psi}^{\text{LxF}, +}$. And the Struijs' ``limiter" is obtained in the following:
\begin{eqnarray}
B^{-} = \frac{\max({\bf \Psi}^{\text{LxF},-}/{\bf \Psi}, 0)}{\sum\limits_{*\in\left\{-, +\right\}} \max({\bf \Psi}^{\text{LxF}, *}/{\bf \Psi},0)},\label{Def:Struijs1d_sys-}\\
B^{+} = \frac{\max({\bf \Psi}^{\text{LxF},+}/{\bf \Psi}, 0)}{\sum\limits_{*\in\left\{-, +\right\}} \max({\bf \Psi}^{\text{LxF}, *}/{\bf \Psi},0)}\label{Def:Struijs_sys+}.
\end{eqnarray}
Let us mention that we calculate $B^{-}$ and $B^{+}$ component by component. Then, we project the ``limiters" $B^{-}$ and $B^{+}$ back to the physical space
\begin{equation}\label{Eq:Struijs_sys1d}
\pmb{\beta}^{-} = \bar{R}B^{-}, ~~ \pmb {\beta}^{+} = \bar{R}B^{+}.
\end{equation}
As for one dimensional systems dissipation residuals, are given according to \eqref{Def:Streamlinediss} as follows:
\begin{equation}
{\bf \Phi}^{-}_{\text{diss}} = -\frac{1}{2}\bar{R}\frac{\bar{\Lambda}}{\left|\bar{\Lambda}\right|}\bar{L}{\bf \Phi}, ~~{\bf \Phi}^{+}_{\text{diss}} = \frac{1}{2}\bar{R}\frac{\bar{\Lambda}}{\left|\bar{\Lambda}\right|}\bar{L}{\bf \Phi},
\end{equation}
where $\bar{\Lambda}$ is the diagonal matrix $\Lambda$ evaluated at the average state and $\left|\bar{\Lambda}\right|$ is the diagonal matrix of the absolutes of all elements in $\bar{\Lambda}$. As before, in order to prevent any element in $\left|\bar{\Lambda}\right|$ from becoming zero, we also take the Roe's correction, as defined in \eqref{Def:Roecorr}.

Hence, we get same formulation as \eqref{Def:rd1d_scl} to distribute the total residual within each cell.
Finally, as in the scalar case, the point value ${\bf u}_{i}$ can be updated in the pseudo time-marching semi-discrete scheme \eqref{Eq:RDscheme}, which is again discretized by a third order TVD Runge-Kutta scheme in our numerical experiments until the steady state is reached.

We now summarize the procedure of the high order RD finite difference WENO-ZQ scheme for one-dimensional steady state systems:
\vspace{0.2cm}
\\
1. Compute the total residual component by component defined in \Cref{Def:residual} using WENO-ZQ integration with a proper accuracy for the source term.
\\
2. 
Project the residuals obtained by the first order Lax-Friedrichs distribution to local characteristic fields, and then obtain the Struijs' ``limiters" \eqref{Def:Struijs1d_sys-}, \eqref{Def:Struijs_sys+}, then project the ``limiters" back to the physical space as in \eqref{Eq:Struijs_sys1d}.
\\
3. Compute the streamline dissipation residuals, then distribute the total residual within each cell according to the SUPG-like principle, which is defined in \eqref{Def:rd1d_scl}. 
\\
4. Update the point values though sending the residuals in the physical space and forward in pseudo time \eqref{Eq:RDscheme} by a TVD Runge-Kutta time discretization until the steady state is reached.

\section{High order RD finite difference WENO-ZQ schemes in two dimension}

In this section, we design a high order RD finite difference WENO-ZQ scheme for two-dimensional steady state problems. We will use Cartesian meshes as examples to describe our algorithm. In the \Cref{Subsec:scalar2d}, we define the total residual within each cell from the integral form, as in \Cref{Def:residual}, and then describe the distribution mechanism. In the \Cref{Subsec:sys2d}, we extend the scheme naturally to two-dimensional systems, based on a local characteristic field decomposition. 

\subsection{Two-dimensional scalar problems}\label{Subsec:scalar2d}

We consider the two-dimensional scalar steady state problem
\begin{equation}\label{Eq:ssp2d}
f(u)_{x}  + g(u)_{y} = s(u, x, y).
\end{equation}
We define the grid to be $\left\{(x_{i}, y_{j})\right\}$, the grid function $u_{ij}$, the cell $I_{i+\frac{1}{2}, j+\frac{1}{2}}=\left[x_{i}, x_{i+1}\right]\times \left[y_{j}, y_{j+1}\right]$, the step x-direction $\Delta x_{i+\frac{1}{2}}$, the step y-direction $\Delta y_{j+\frac{1}{2}}$, the control volume centered at $(x_{i}, y_{j})$ to be $C_{ij}$(formed by connecting the centers of the four cells sharing $(x_{i},y_{j})$ as a common node), and the area of $C_{ij}$ is denoted by $\left|C_{ij}\right|$.

The total residual in the cell $I_{i+\frac{1}{2}, j+\frac{1}{2}}$ is defined by
\begin{equation}\label{Def:residual2d}
\begin{array} {lcl} 
 \Phi_{i+\frac{1}{2},j+\frac{1}{2}} & = & \int^{y_{j+1}}_{y_{j}}\int^{x_{i+1}}_{x_{i}}\!(f(u)_{x} + g(u)_{y} - s(u, x, y))\,\mathrm{d}x\,\mathrm{d}y \\
& = & \int^{y_{j+1}}_{y_{j}}\!(f(u(x_{i+1}, y)) - f(u(x_{i}, y)))\,\mathrm{d}y + \int^{x_{i+1}}_{x_{i}}\!(g(u(x, y_{j+1})) - g(u(x, y_{j})))\,\mathrm{d}x \\
& & -\int^{y_{j+1}}_{y_{j}}\int^{x_{i+1}}_{x_{i}}\!s(u(x, y), x, y)\,\mathrm{d}x\,\mathrm{d}y.
\end{array}
\end{equation}
If we can reach the zero residual limit, i.e., if $\Phi_{i+\frac{1}{2}, j+\frac{1}{2}} = 0$ for all i and j, the accuracy of the scheme is determined by the accuracy of the approximation to the integrations of the fluxes and the source term.

To approximate the integrations of the fluxes, which are one-dimensional integrals, we use a fourth order WENO-ZQ integration described in the \Cref{Subsec:scalar1d}. As for the source term $\int^{y_{j+1}}_{y_{j}}\int^{x_{i+1}}_{x_{i}}\!s(u, x, y)\, \mathrm{d}x\, \mathrm{d}y$, we can approximate it in a dimension by dimension fashion, which is explained as follows: 

First, we define
$$S_{j+\frac{1}{2}}(x) = \int^{y_{i+1}}_{y_{i}}\!s(u(x, y), x, y)\,\mathrm{d}y,$$
and then
$$\int^{y_{j+1}}_{y_{j}}\int^{x_{i+1}}_{x_{i}}\!s(u, x, y)\,\mathrm{d}x\,\mathrm{d}y = \int^{x_{i+1}}_{x_{i}}\!S_{j+\frac{1}{2}}(x)\,\mathrm{d}x.$$
The integral $\int^{x_{i+1}}_{x_{i}}\!S_{j+\frac{1}{2}}(x)\,\mathrm{d}x$ can be approximated by a fourth order WENO-ZQ integration in the x-direction, using $\left\{S_{j+\frac{1}{2}}(x_{i+k})\right\}_{k=-1,\cdots,2}$. By the definition of $S_{j+\frac{1}{2}}(x)$, $S_{j+\frac{1}{2}}(x_{i+k})$ can again be approximated by a fourth order WENO-ZQ integration in the y-direction, using $\left\{s(u_{i+k,j+l}, x_{i+k}, y_{j+l})\right\}_{l=-1,\cdots,2}$. Therefore, the integration of the source term can be approximated dimension by dimension, and the fourth order accuracy is the zero residual limit. 

Next, we start to distribute the total residuals. In the cell $I_{i+\frac{1}{2}, j+\frac{1}{2}}=\left[x_{i},x_{i+1}\right]\times\left[y_{j},y_{j+1}\right]$, the total residual is $\Phi_{i+\frac{1}{2}, j+\frac{1}{2}}$, and it is to be distributed to the vertices of the cell, which are defined to be $M_{1} = (x_{i+1}, y_{j+1})$, $M_{2} = (x_{i+1}, y_{j})$, $M_{3} = (x_{i}, y_{j+1})$ and $M_{4} = (x_{i}, y_{j})$. Here we denote the residuals distributed to the vertices $M_{k}$ as $\Phi^{k}_{i+\frac{1}{2}, j+\frac{1}{2}}$, $k=1,2,3,4$. For simplicity and without ambiguity, we drop the subscript $(i+\frac{1}{2}, j+\frac{1}{2})$ off in the notations. For the conservation and the residual property, we require $\Phi = \sum^{4}_{k=1}\Phi^{k}$ and $\left|\Phi^{k}\right|/\left|\Phi\right|$ to be uniformly bounded.

To have a SUPG-like scheme, one way to distribute the total residual $\Phi$ is the following:

Step 1. First order Lax-Friedrichs linear distribution is given by
\begin{equation}\label{Def:LaxF2dscl}
\begin{split}
\Phi^{\text{LxF}, M_{1}} & = \frac{1}{4}\Phi+\alpha(u_{i+1, j+1} - \bar{u}),\\
\Phi^{\text{LxF}, M_{2}} & = \frac{1}{4}\Phi+\alpha(u_{i+1, j} - \bar{u}),\\
\Phi^{\text{LxF}, M_{3}} & = \frac{1}{4}\Phi+\alpha(u_{i, j+1} - \bar{u}),\\
\Phi^{\text{LxF}, M_{4}} & = \frac{1}{4}\Phi+\alpha(u_{i, j} - \bar{u}),
\end{split}
\end{equation}
where $\bar{u}$ is an average state in the cell taken to be $\frac{1}{4}(u_{i+1, j+1} + u_{i+1, j} + u_{i, j+1} + u_{i, j})$, and $\alpha$ is determined by
\begin{equation}
\alpha = \max (\Delta x_{i+\frac{1}{2}}, \Delta y_{j+\frac{1}{2}})\cdot \max\limits_{u_{ij}\in I_{i+\frac{1}{2}, j+\frac{1}{2}}}\left\{\left| f'(u_{ij})\right| + \left|g'(u_{ij})\right|\right\}.
\end{equation}
The Struijs' ``limiter" is given by
\begin{equation}\label{Def:Struijs2d}
\beta^{M_{k}} = \frac{\max(\Phi^{\text{LxF}, M_{k}}/\Phi, 0)}{\sum\limits_{M{*}\in I_{i+\frac{1}{2}, j+\frac{1}{2}}}\max(\Phi^{\text{LxF}, M_{*}}/\Phi, 0)}, ~~ k=1, \cdots, 4.
\end{equation}

Step 2. The streamline dissipation term is defined by
\begin{equation}\label{Def:Streamlinediss2d}
\int_{I_{i+\frac{1}{2}, j+\frac{1}{2}}}\!(\nabla_{u}f(u)\cdot \nabla\varphi^
{M_{k}})\tau (\nabla_{u}f(u)\cdot\nabla u - s(u, x, y))\, \mathrm{d}x\, \mathrm{d}y,
\end{equation}
where $\varphi^{M_{k}}$ is the basis function associated to the node $M_{k}$, $k = 1, \cdots, 4$  in the cell $I_{i+\frac{1}{2}, j+\frac{1}{2}}$, and we take them as follows:
\begin{equation}
\begin{split}
\varphi^{M_{1}} & = \frac{x - x_{i}}{x_{i+1} - x_{i}}\frac{y-y_{j}}{y_{j+1} - y_{j}}, \\
\varphi^{M_{2}} & = \frac{x - x_{i}}{x_{i+1} - x_{i}}\left(1 - \frac{y-y_{j}}{y_{j+1} - y_{j}}\right), \\
\varphi^{M_{3}} & = \frac{y-y_{j}}{y_{j+1} - y_{j}} \left(1- \frac{x - x_{i}}{x_{i+1} - x_{i}}\right),\\
\varphi^{M_{4}} & = \left(1 - \frac{x-x_{i}}{x_{i+1} - x_{i}}\right)\left(1 - \frac{y-y_{j}}{y_{j+1} - y_{j}}\right).
\end{split}
\end{equation}
And $\tau^{-1}$ is taken to be
\begin{equation}
\tau^{-1} = \sum\limits_{M_{k}\in I_{i+\frac{1}{2},j+\frac{1}{2}}}\left|(f'(\bar{u}), g'(\bar{u}))\cdot\nabla\varphi^{M_{k}}(x^{M_{k}}, y^{M_{k}})\right|.
\end{equation}
As for two dimensional scalar case, we take 
\begin{equation}\label{Def:Dissipation2d}
\Phi^{k}_{\text{diss}} = (f'(\bar{u}), g'(\bar{u}))\cdot \nabla\varphi^{M_{k}}(x^{M_{k}}, y^{M_{k}})\tau \Phi, ~~k = 1, \cdots, 4.
\end{equation}
Hence, we get the way to distribute the total residual within each cell as follows:
\begin{equation}\label{Def:rd2d_scl}
\Phi^{k} = \beta^{M_{k}}\Phi + \Phi^{k}_{\text{diss}}, ~~k = 1, \cdots, 4.
\end{equation}

The point value $u_{ij}$ is then updated through sending the distributed residuals to the point $(x_{i}, y_{j})$, as in a pseudo time-marching scheme, which can be written as a semi-discrete system
\begin{equation}\label{Eq:RDscheme2d}
\frac{\mathrm{d}u_{ij}}{\mathrm{d}t} + \frac{1}{\left|C_{ij}\right|}\left(\Phi^{1}_{i-\frac{1}{2}, j-\frac{1}{2}} + \Phi^{2}_{i-\frac{1}{2}, j+\frac{1}{2}} + \Phi^{3}_{i+\frac{1}{2}, j-\frac{1}{2}} + \Phi^{4}_{i+\frac{1}{2}, j+\frac{1}{2}}\right) = 0.
\end{equation}
We again use a third order TVD Runge-Kutta scheme for the pseudo time discretization, as in the one-dimensional case.

We now summarize the procedure of the high order RD finite difference WENO-ZQ scheme for two dimensional scalar steady state problems:
\vspace{0.2cm}
\\
1. Compute the total residual defined in \Cref{Def:residual2d} using WENO-ZQ intergration dimension by dimension with a proper accuracy for the source term.\\
2. Distribute the total residual within each cell according to the SUPG-like principle, which is defined in \eqref{Def:LaxF2dscl}, \eqref{Def:Struijs2d}, \eqref{Def:Dissipation2d}, \eqref{Def:rd2d_scl}.\\
3. Update the point values through sending the residuals and forward in pseudo time \eqref{Eq:RDscheme2d} by a TVD Runge-Kutta time discretization until the steady state is reached.

\subsection{Two-dimensional systems}\label{Subsec:sys2d}

Consider a two-dimensional steady state system \eqref{Eq:ssp2d} where {\bf u}, {\bf f}({\bf u}), {\bf g}({\bf u}) and {\bf s}({\bf u}, x, y) are vector-valued functions in $\mathbb{R}^{m}$. For hyperbolic systems, we assume that any real linear combination of the Jacobians $n_{x}{\bf f}'({\bf u}) + n_{y}{\bf g}'({\bf u})$ is diagonalizable with real eigenvalues. In particular, we assume ${\bf f}'({\bf u})$ and ${\bf g}'({\bf u})$ can be written as $R_{x}\Lambda_{x}L_{x}$ and $R_{y}\Lambda_{y}L_{y}$, respectively, where $\Lambda_{x}$ and $\Lambda_{y}$ are diagonal matrices with real eigenvalues on the diagonal, and $L_{x}$, $R_{x}$ and $L_{y}$, $R_{y}$ are matrices of left and right eigenvectors for the corresponding Jacobians.

The grid, grid function, step x-direction, step y-direction and control volume are denoted as in the \Cref{Subsec:scalar2d}. The total residual in the cell $I_{i+\frac{1}{2}, j+\frac{1}{2}} = \left[x_{i}, x_{i+1}\right]\times\left[y_{j}, y_{j+1}\right]$ is still defined by \eqref{Def:residual2d}. As before, if we can reach the zero residual limit of the scheme, the accuracy of the scheme is determined by the accuracy of the approximations to the integrations of the fluxes and the source term. We again use a fourth order WENO-ZQ integration described in the \Cref{Subsec:scalar1d}. For simplicity and without ambiguity, we drop the subscript $(i+\frac{1}{2}, j+\frac{1}{2})$ off in the notations in the following.

We distribute the total residual ${\bf \Phi}$ to the four vertices $\left\{M_{k}\right\}_{k=1,\cdots, 4}$, which is defined in the \Cref{Subsec:scalar2d} and the corresponding residuals are still denoted by $\left\{{\bf \Phi}^{k}\right\}_{k=1,\cdots, 4}$, where ${\bf \Phi}^{k}\in \mathbb{R}^{m}$. We also require ${\bf \Phi} = \sum^{4}_{k=1}{\bf \Phi}^{k}$ and the residual property that $\left|{\bf \Phi}^{k}\right|/\left|{\bf \Phi}\right|$ should be uniformly bounded. First, we compute an average state $\bar{{\bf u}}$ in $I_{i+\frac{1}{2}, j+\frac{1}{2}}$, using either arithmetic mean or Roe's average \citep{Roe_JCP1981}. And then denote $\bar{L}$ and $\bar{R}$ as the matrices with left and right eigenvectors $L$ and $R$ of $n_{x}{\bf f}'({\bf u}) + n_{y}{\bf g}'({\bf u})$ evaluated at the average state, where ${\bf n}=(n_{x}, n_{y})$ can be any direction. The first order Lax-Friedrichs linear distribution is again defined by \eqref{Def:LaxF2dscl}. Then we project ${\bf \Phi}^{\text{LxF}, M_{k}}, k=1,\cdots, 4$ to the characteristic fields, namely, ${\bf \Psi}^{\text{LxF}, M_{k}} = \bar{L}{\bf \Phi}^{\text{LxF}, M_{k}}, k=1, \cdots, 4$, with ${\bf \Psi} = \sum^{4}_{k=1}{\bf \Psi}^{\text{LxF}, M_{k}}$. The Struijs' ``limiter" is obtained in the following:
\begin{equation}
B^{M_{k}} = \frac{\max({\bf \Psi}^{\text{LxF},M_{k}}/{\bf \Psi}, 0)}{\sum\limits_{M_{*}\in I_{i+\frac{1}{2}, j+\frac{1}{2}}} \max({\bf \Psi}^{\text{LxF}, M_{*}}/{\bf \Psi},0)}, ~~ k = 1,\cdots, 4.
\end{equation}
Let us mention that we compute $B^{M_{k}}$ component by component.
Then, we project the ``limiters" $B^{M_{k}}, k= 1,\cdots, 4$ back to the physical space, we obtain
\begin{equation}\label{Eq:Struijs_sys2d}
\pmb {\beta}^{M_{k}} = \bar{R}B^{M_{k}}, ~~ k = 1,\cdots, 4.
\end{equation}
As for two dimensional systems streamline dissipation residuals, are given according to \eqref{Def:Streamlinediss2d} as follows:
\begin{equation}\label{Def:Streamlinediss2dsys}
{\bf \Phi}^{k}_{\text{diss}}= ({\bf f}'(\bar{{\bf u}}), {\bf g}'(\bar{{\bf u}})
)\cdot\nabla\varphi^{M_{k}}(x^{M_{k}}, y^{M_{k}})\tau {\bf \Phi},~~ k = 1,\cdots,4.
\end{equation}
Here $\tau^{-1}$ is taken to be
\begin{equation}
\tau^{-1} = \sum\limits_{M_{k}\in I_{i+\frac{1}{2},j+\frac{1}{2}}}\left|({\bf f}'(\bar{{\bf u}}), {\bf g}'(\bar{{\bf u}}))\cdot\nabla\varphi^{M_{k}}(x^{M_{k}}, y^{M_{k}})\right|.
\end{equation}

Hence, we get same formulation as \eqref{Def:rd2d_scl} to distribute the total residual within each cell. Finally, as in the two-dimensional scalar case, the point value ${\bf u}_{ij}$ can be updated in the pseudo time-marching semi-discrete scheme \eqref{Eq:RDscheme2d}, which is again discretized by a third order TVD Runge-Kutta scheme in our numerical experiments until the steady state is reached.

We now summarize the procedure of the high order RD finite difference WENO-ZQ scheme for two-dimensional steady state systems:
\vspace{0.2cm}\\
1. Compute the total residual component by component defined in \Cref{Def:residual2d} using WENO-ZQ integration dimension by dimension with a proper accuracy for the source term.
\\
2. Project the residuals obtained by the first order Lax-Friedrichs distribution to local characteristic fields, then obtain the Struijs' ``limiters", then project the ``limiters" back to the physical space as in \eqref{Eq:Struijs_sys2d}.
\\
3. Compute dissipation residuals, then distribute the total residual within each cell according to the SUPG-like principle, which is  defined in \eqref{Def:Streamlinediss2dsys}.
\\
4. Update the point values though sending the residual in the physical space and forward in pseudo time \eqref{Eq:RDscheme2d} by a TVD Runge-Kutta time discretization until the steady state is reached.

\section{Numerical results}

In this section, we present the numerical results of the proposed fourth order residual distribution finite difference WENO-ZQ method for hyperbolic conservation laws with source terms in scalar and system test problems in one and two dimensions. Pseudo time discretization  towards steady state 
is by the third order TVD Runge-Kutta method in  all numerical simulations.

All the spatial discretizations in our numerical results are uniform. And the CFL numbers are taken to be 0.3 in all problems. We remark here that the choice of the CFL number certainly affects the number of iteration to reach a steady state, but here we choose it to be sufficiently large while maintaining stability for all cases.

\subsection{The one-dimensional scalar problems}

In this section, all numerical steady state is obtained with $L^{1}$ residue reduced to the round-off level.
\begin{exmp}\label{Eg:bur1d_sm}
We solve the steady state solution of the one-dimensional Burgers equation with a source term:

\begin{equation}
u_{t}+\left(\frac{u^{2}}{2}\right)_{x}=\sin x \cos x
\end{equation}
with the initial condition
\begin{equation}\label{Eq:bur1d_sm_init}
u(x,0) = \beta \sin x
\end{equation}
and the boundary condition $u(0, t)=u(\pi, t) = 0$. This problem was studied in \cite{S.A.Gottlieb_ANM1986} as an example of multiple steady state solutions for characteristic initial value problems. The steady state solution to this problem depends on the value of $\beta$: if $-1<\beta< 1$, a shock will form within the domain $\left[0,\pi\right]$; otherwise, the solution will be smooth at first, followed by a shock forming at the boundary $x = \pi ~(\beta \geq 1)$ or $x = 0~(\beta \leq -1)$, and later converge to a smooth steady state $u(x, \infty) = \sin x~(\beta \geq 1)$ or $u(x, \infty) = -\sin x~(\beta \leq -1)$, respectively. In order to test the order of accuracy, we take $\beta =2$ to have a smooth stationary solution.
The numerical results are shown in the \Cref{Tab:bur1d_sm}. We can see clearly that the fourth order is reached on the uniform meshes.

\begin{table}[ht!]\small
\caption{Errors and numerical orders of accuracy for the fourth order SUPG-like RD finite difference WENO-ZQ scheme for the \Cref{Eg:bur1d_sm} on uniform meshes with N cells}
\vspace{0.1 in}
\centering
\begin{tabular}{lllll}
\hline
\text{$N$} & \text{$L^{1}$ error} & \text{Order} & \text{$L^{\infty}$ error} & \text{Order} \\
\hline
      20 &     3.96E-05 & &     6.45E-05 & \\ 
     40 &     2.77E-06 &     3.84 &     4.49E-06 &     3.84 \\ 
     80 &     1.81E-07 &     3.94 &     2.88E-07 &     3.96 \\ 
    160 &     1.15E-08 &     3.98 &     1.81E-08 &     3.99 \\ 
    320 &     7.21E-10 &     3.99 &     1.13E-09 &     4.00 \\ 
    640 &     4.52E-11 &     4.00 &     7.10E-11 &     4.00 \\ \hline
\end{tabular}
\label{Tab:bur1d_sm}
\end{table}
\end{exmp}

\begin{exmp}\label{Eg:bur1d_nsm}
We consider the same problem as the \Cref{Eg:bur1d_sm}, but take $\beta = 0.5$ in the initial condition \eqref{Eq:bur1d_sm_init}. As mentioned in the previous example, when $-1 < \beta < 1$, a shock will form within the domain, which separates two branches ($\sin x$ and $-\sin x$) of the steady state. The location of the shock is determined by the parameter $\beta$ through conservation of mass ($\int^{\pi}_{0} u \,dx = 2\beta$), and can be derived to be $\pi - \arcsin \sqrt{1-\beta^{2}}$. For the case $\beta = 0.5$, the shock location is approximately 2.0944. The numerical solution on the uniform meshes is shown in the \Cref{Fig:bur1d_nsm}. We can see that the numerical shock is at the correct location and is resolved well.
{\color{red} We also observe the convergence histories by different CFL numbers and the results are shown in \Cref{Fig:bur1d_nsm_conv}. We can see that the CFL number influences the convergence history, the larger CFL number and the faster convergence. When CFL number is 0.7, the $L^{1}$ residue stagnates only at $10^{-7}$ level. }

\begin{figure}[ht!]
\centering
\includegraphics[width=7cm]{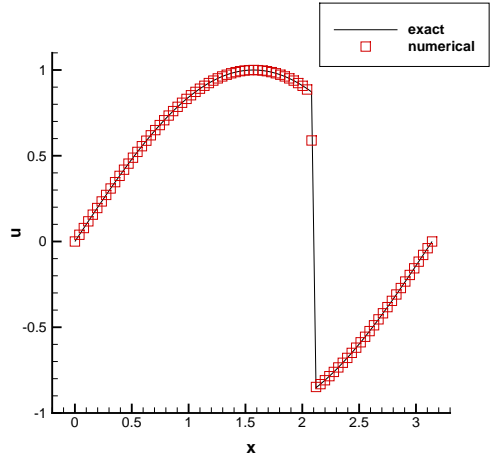}
\caption{The numerical solution (symbols) versus the exact solution (solid line) for the \Cref{Eg:bur1d_nsm} with 80 cells}\label{Fig:bur1d_nsm}
\end{figure}

\begin{figure}[ht!]
\centering
\begin{minipage}[b]{0.4\textwidth}
\centering
\includegraphics[width=6.5cm]{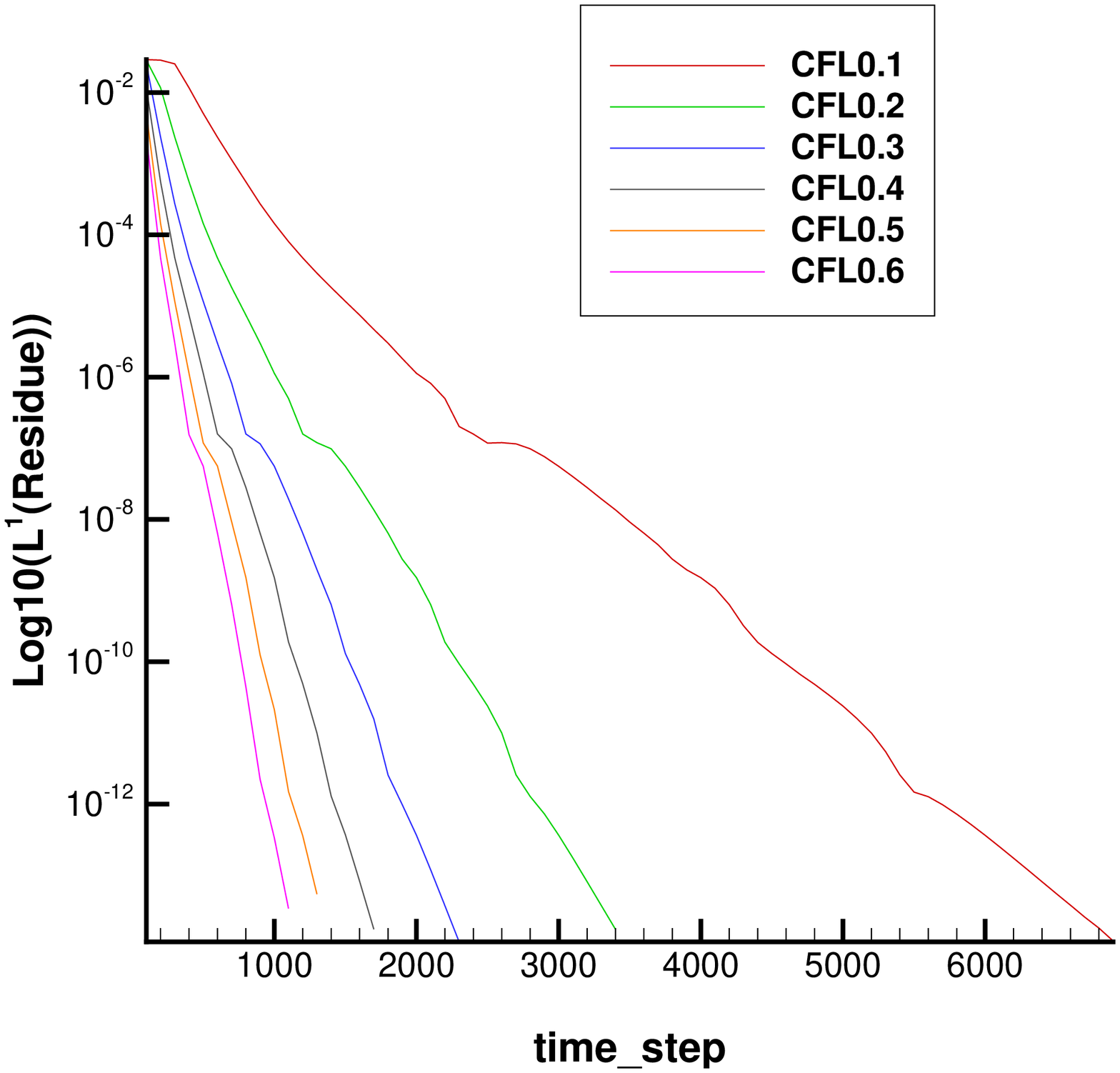}\\
\end{minipage}
\hspace{0.4cm}
\begin{minipage}[b]{0.4\textwidth}
\centering
\includegraphics[width=6.5cm]{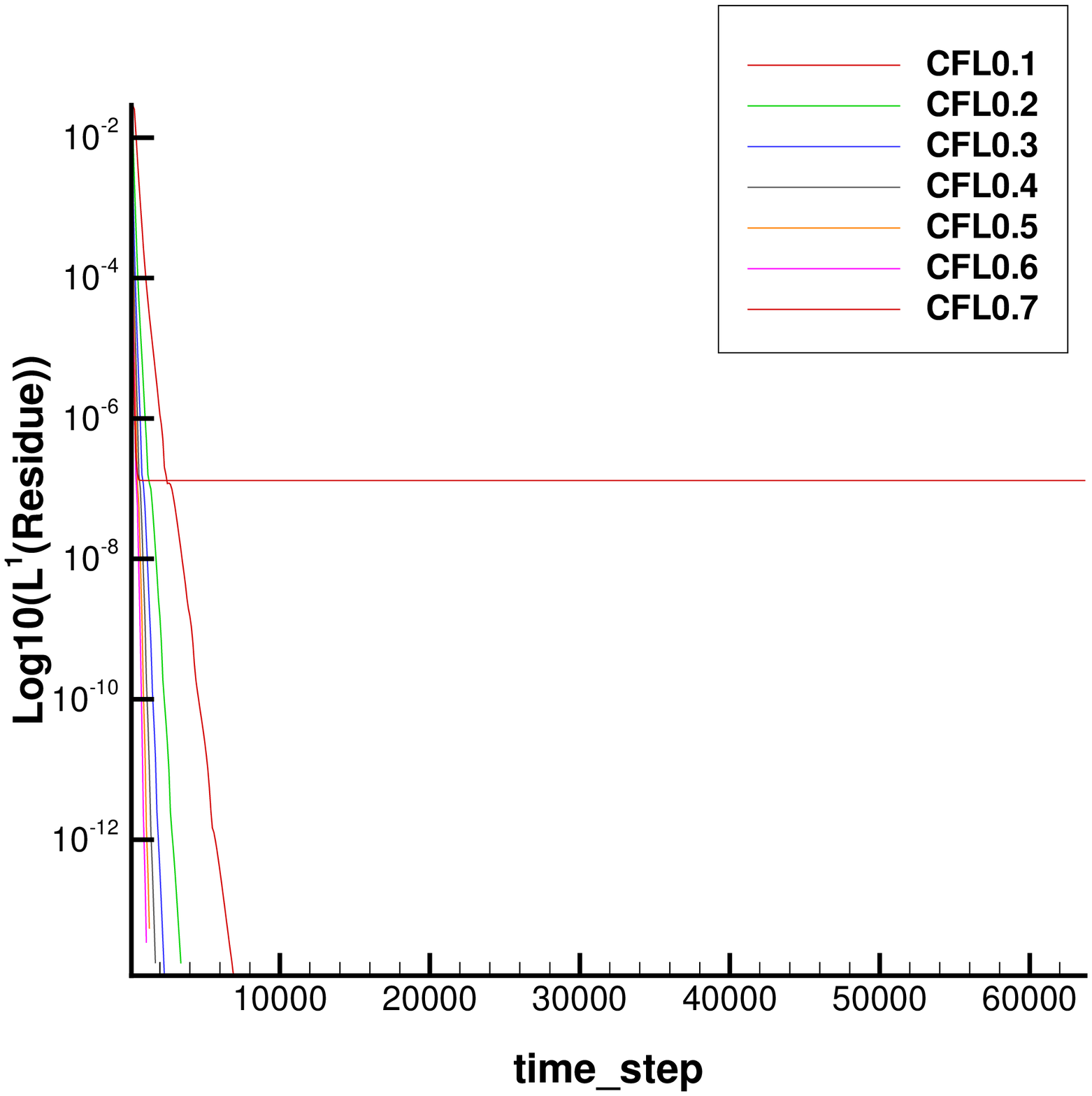}\\
\end{minipage}

\caption{The convergence histories of $L^{1}$ residue for the \Cref{Eg:bur1d_nsm}}\label{Fig:bur1d_nsm_conv}
\end{figure}
\end{exmp}

\begin{exmp}\label{Eg:bur1d_nsm2}
We consider the steady state solutions of the Burgers equation with a different source term, which depends on the solution itself:
\begin{equation}
u_{t} + \left( \frac{u^{2}}{2}\right)_{x} = -\pi \cos(\pi x)u, ~~ x\in \left[0, 1\right]
\end{equation}
equipped with the boundary conditions $u(0, t) = 1$ and $u(1, t) = -0.1$. This problem has two steady state solutions with shocks
\[
u(x, \infty) = \begin{cases}
u^{+} = 1 - \sin(\pi x) & \text{if $0 \leq x < x_{s}$}, \\
u^{-} = -0.1 - \sin(\pi x) & \text{if $x_{s} \leq x \leq 1$},
\end{cases}
\]
where $x_{s} = 0.1486$ or $x_{s} = 0.8514$. Both solutions satisfy the Rankine-Hugoniot jump condition and the entropy conditions, but only the one with the shock at 0.1486 is stable for a small perturbation. This problem was studied in \citep{Embid.Goodman.Majda_SIAM_SSC1984} as an example of multiple steady states for one-dimensional transonic flows. This case is tested to demonstrate that starting with a reasonable perturbation of the stable steady state, the numerical solution converges to the stable one.

The initial condition is given by
\[
u(x, 0) = \begin{cases}
\hfil 1 & \text{if $0 \leq x < 0.5$},\\
\hfil -0.1  & \text{if $0.5 \leq x \leq 1$},
\end{cases}
\]
where the initial jump is located in the middle of the position of the shocks in the two admissible steady state solution. The numerical result and the exact solution are displayed in the \Cref{Fig:bur1d_nsm2}. We can see the correct shock location and good resolution of the shock.
{\color{red}We also observe the convergence histories by different CFL numbers and the results are shown in \Cref{Fig:bur1d_nsm2_conv}. We can see that the CFL number influences the convergence history, the larger CFL number and the faster convergence. When CFL number is 0.7, the $L^{1}$ residue stagnates only at $10^{-12}$ level.}

\begin{figure}[ht!]
\centering
\includegraphics[width=7cm]{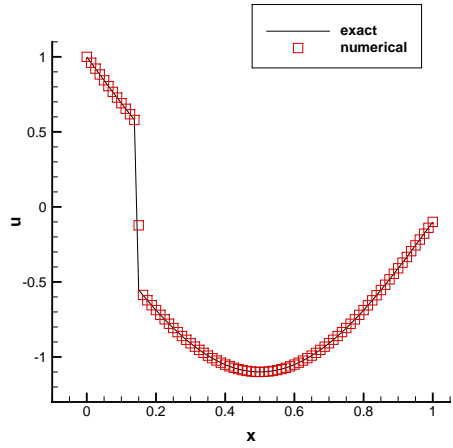}
\caption{The numerical solution (symbols) versus the exact solution (solid line) for the \Cref{Eg:bur1d_nsm2} with 80 cells}\label{Fig:bur1d_nsm2}
\end{figure}

\begin{figure}[ht!]
\centering
\begin{minipage}[b]{0.4\textwidth}
\centering
\includegraphics[width=6.5cm]{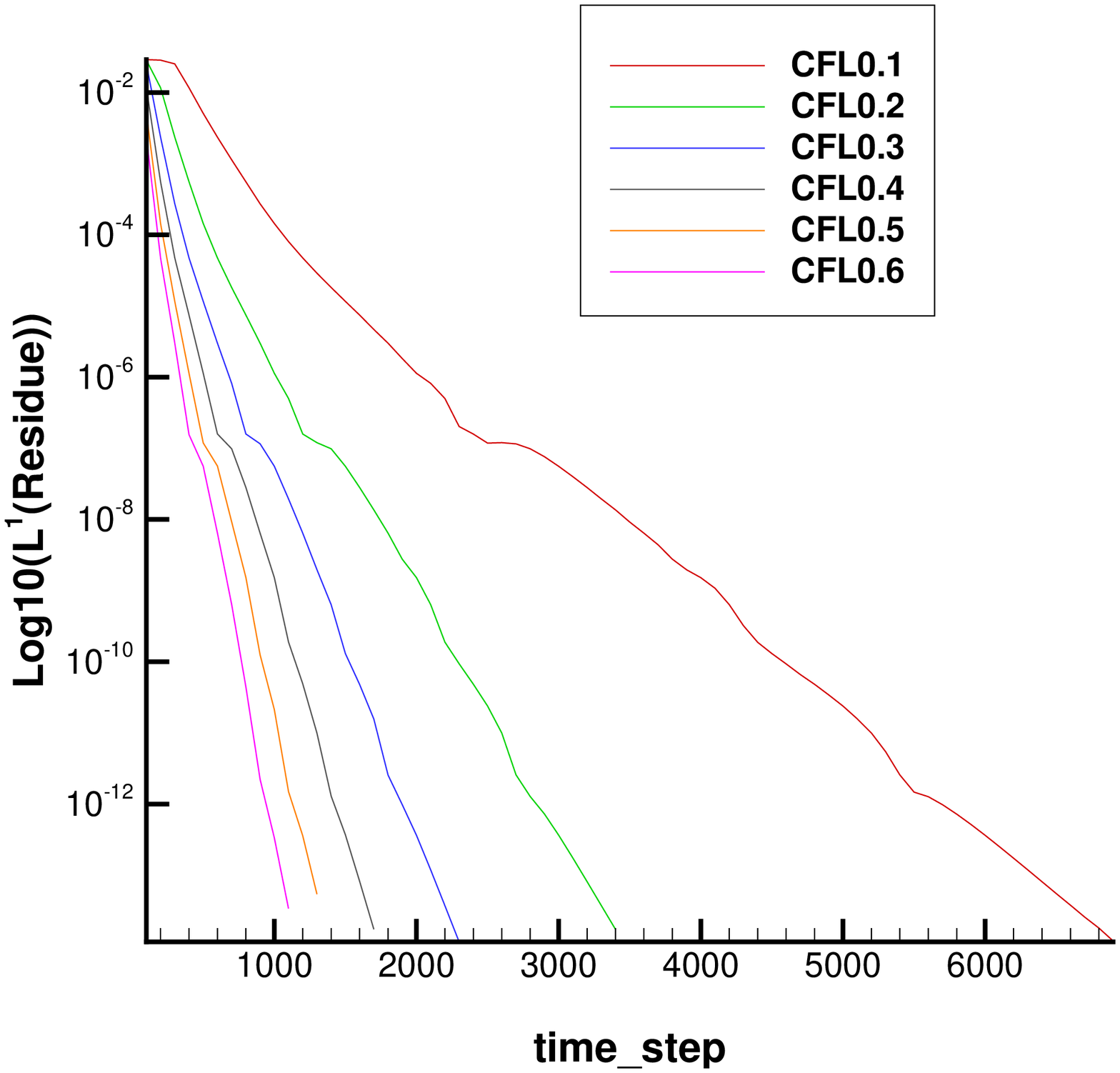}\\
\end{minipage}
\hspace{0.4cm}
\begin{minipage}[b]{0.4\textwidth}
\centering
\includegraphics[width=6.5cm]{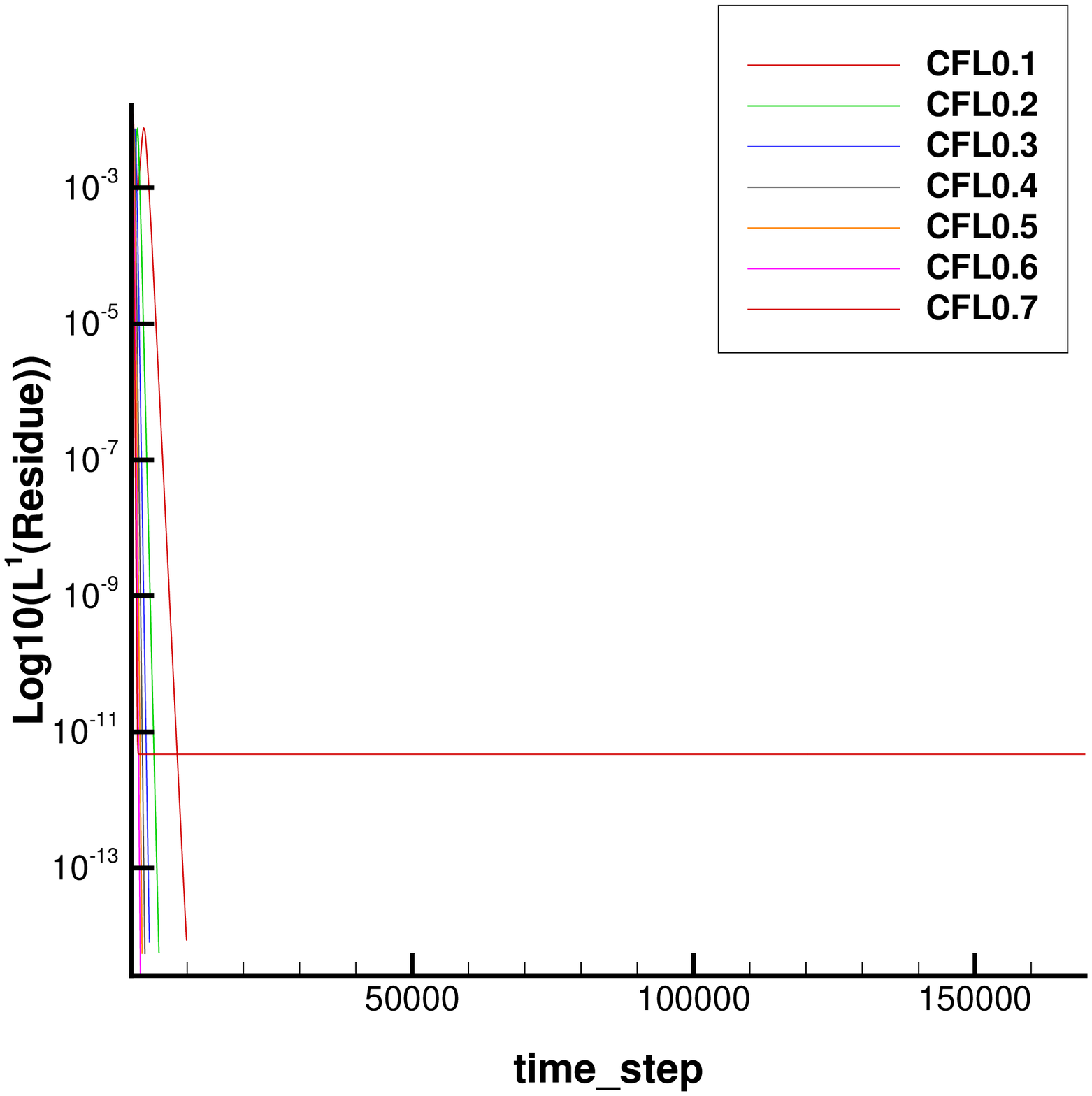}\\
\end{minipage}

\caption{The convergence histories of  $L^{1}$ residue for the \Cref{Eg:bur1d_nsm2}}\label{Fig:bur1d_nsm2_conv}
\end{figure}

\end{exmp}

\subsection{The one-dimensional systems}

\begin{exmp}\label{Eg:shallowwater}
We solve the steady state solutions of the one-dimensional shallow water equation

\begin{equation}
\left(\begin{array}{c}
h \\
hu
\end{array}\right)_{t}
+ \left(\begin{array}{c}
hu \\
hu^{2}+\frac{1}{2}gh^{2}

\end{array}\right)_{x}
=\left(\begin{array}{c}
0 \\
-ghb_{x}
\end{array}\right),
\end{equation}
where h denotes the water height, u is the velocity of the fluid, b(x) represents the bottom topography and g is the gravitational constant.

Starting from a stationary initial condition, which itself is a steady state solution, we can check the order of accuracy. The smooth bottom topography is given by
\begin{equation*}
b(x) = 5\exp^{-\frac{2}{5}(x - 5)^{2}}, ~~x\in [0, 10].
\end{equation*}
The initial condition is the stationary solution
\begin{equation*}
h + b = 10, ~~hu = 0
\end{equation*}
and the exact steady state solution is imposed as the boundary condition.

We test our scheme on uniform meshes. The numerical results are shown in the \Cref{Tab:shallowwater}. We can clearly see the order of accuracy and the errors.

\begin{table}[ht!]\small
\caption{Errors and numerical orders of accuracy for the water height h of the fourth order SUPG-like RD finite difference WENO-ZQ scheme for the \Cref{Eg:shallowwater} on uniform meshes with N cells}
\vspace{0.1 in}
\centering
\begin{tabular}{lllll}
\hline
\text{$N$} & \text{$L^{1}$ error} & \text{Order} & \text{$L^{\infty}$ error} & \text{Order} \\
\hline
     20 &     3.43E-02 & &     1.08E-02 & \\
     40 &     9.02E-03 &     1.93 &     3.28E-03 &     1.72 \\ 
     80 &     2.89E-04 &     4.96 &     1.08E-04 &     4.92 \\ 
    160 &     6.38E-05 &     2.18 &     2.37E-05 &     2.19 \\ 
    320 &     9.04E-07 &     6.14 &     3.29E-07 &     6.17 \\ 
    640 &     7.60E-08 &     3.57 &     2.78E-08 &     3.56 \\ 
   1280 &     1.25E-09 &     5.93 &     4.36E-10 &     6.00 \\ 
   2560 &     7.72E-11 &     4.01 &     2.70E-11 &     4.02 \\ \hline
\end{tabular}
\label{Tab:shallowwater}
\end{table}

\end{exmp}

\begin{exmp}\label{Eg:nozzle1d}
We test our scheme on the steady state solution of the one-dimensional nozzle flow problem
\begin{equation}
\left(\begin{array}{c}
\rho  \\
\rho u \\
E
\end{array}\right)_{t}
+ \left(\begin{array}{c}
\rho u \\
\rho u^{2} + p \\
u(E + p)
\end{array}\right)_{x}
=-\frac{A'(x)}{A(x)}\left(\begin{array}{c}
\rho u \\
\rho^{2}u^{2}/\rho \\
u(E + p)
\end{array}\right), ~~ x\in \left[0, 1\right],
\end{equation}
where $\rho$ denotes the density, u is the velocity of the fluid, E is the total energy, $\gamma$ is the gas constant, which is taken as 1.4, $p = (\gamma -1)(E - \frac{1}{2}\rho u^{2})$ is the pressure, and $\mathit{A(x)}$ represents the area of the cross-section of the nozzle.

We start with an isentropic initial condition, with a shock at $x = 0.5$. The density $\rho$ and pressure $\mathit{p}$ at $-\infty$ are 1, and the inlet Mach number at $x = 0$ is 0.8. The outlet Mach number at $x = 1$ is 1.8, with linear Mach number distribution before and after the shock. The area of the cross-section $\mathit{A(x)}$ is then determined by the relation
\begin{equation*}
A(x)f(\text{Mach number at x}) = \text{constant}, ~~ \forall x\in [0, 1],
\end{equation*}
where
\begin{equation*}
f(w) = \frac{w}{(1+\delta w^{2})^{p}},~~\delta=\frac{1}{2}(\gamma -1),~~p = \frac{1}{2}\cdot\frac{\gamma + 1}{\gamma - 1}.
\end{equation*}

From the \Cref{Fig:nozzle1d}, we can clearly see that the shock is resolved well. {\color{red}We also observe the convergence histories by different CFL numbers and the results are shown in \Cref{Fig:nozzle_conv}. We can see that the CFL number influences the convergence history, the larger CFL number and the faster convergence. When CFL number is 0.8 or 0.9, the $L^{1}$ residue stagnates only at $10^{-6}$ level.}

\begin{figure}[ht!]
\centering
\begin{minipage}[b]{0.4\textwidth}
\centering
\includegraphics[width=6.5cm]{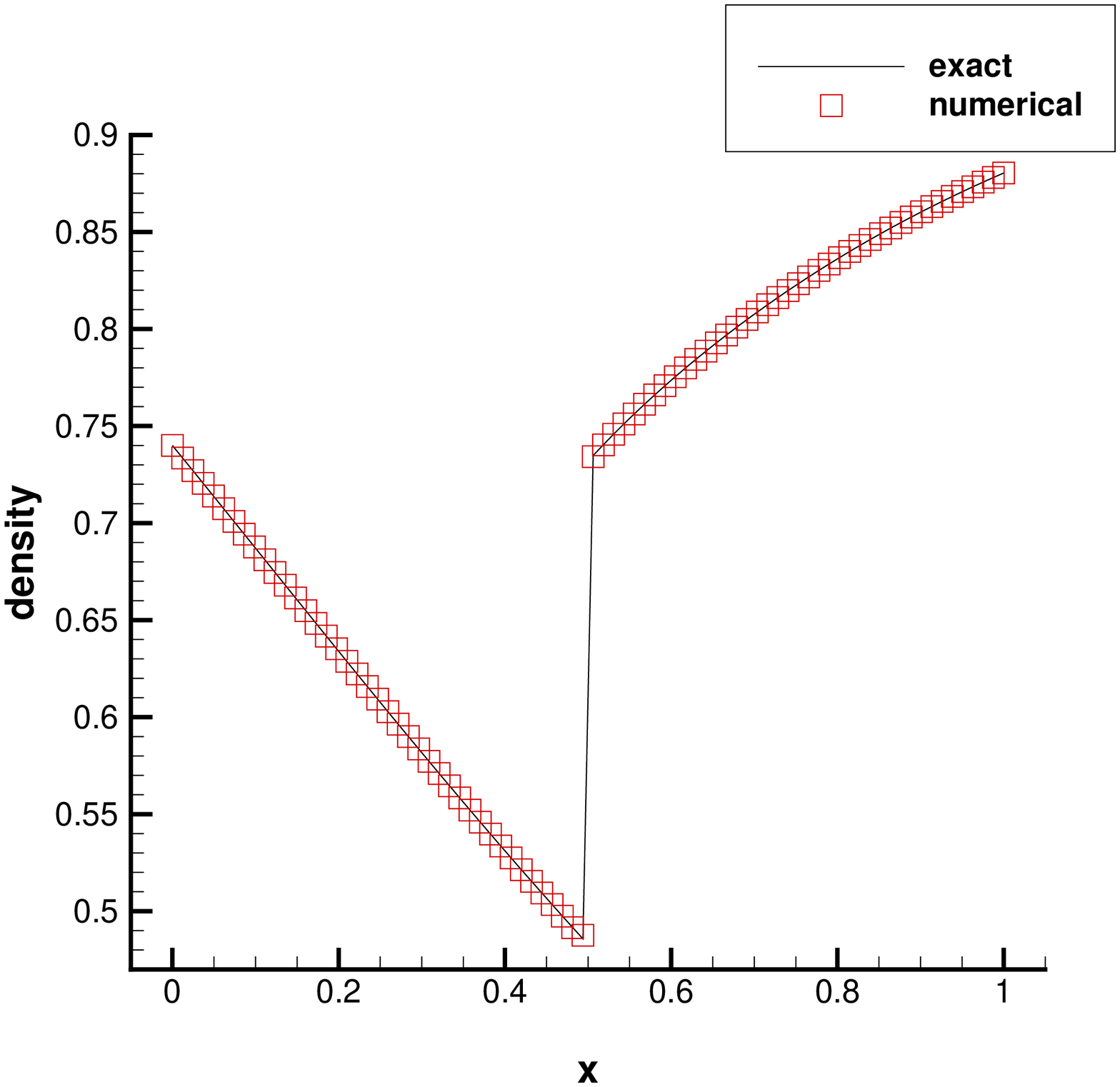}\\
\end{minipage}
\hspace{0.4cm}
\begin{minipage}[b]{0.4\textwidth}
\centering
\includegraphics[width=6.5cm]{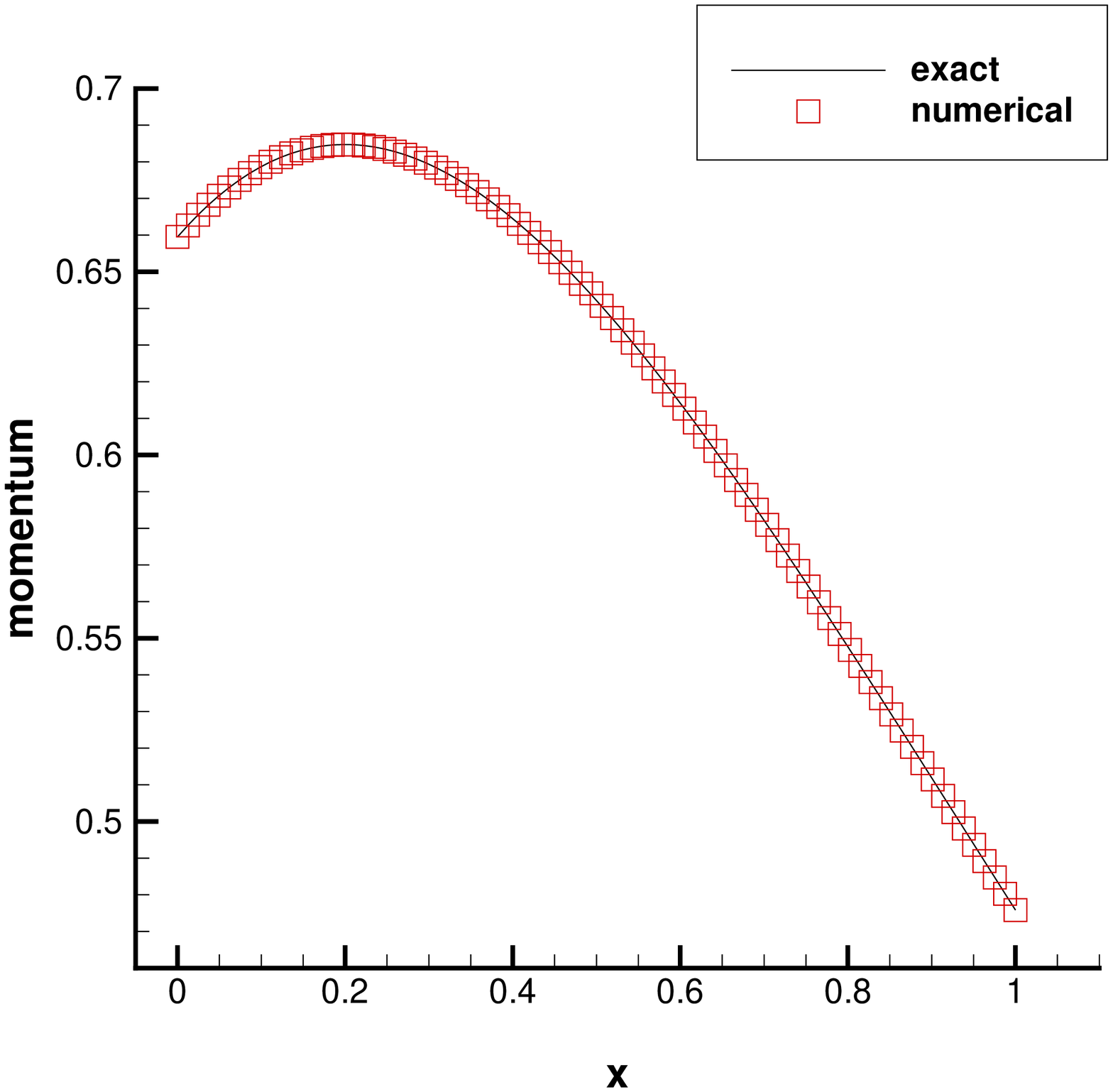}\\
\end{minipage}

\begin{minipage}[b]{0.4\textwidth}
\centering
\includegraphics[width=6.5cm]{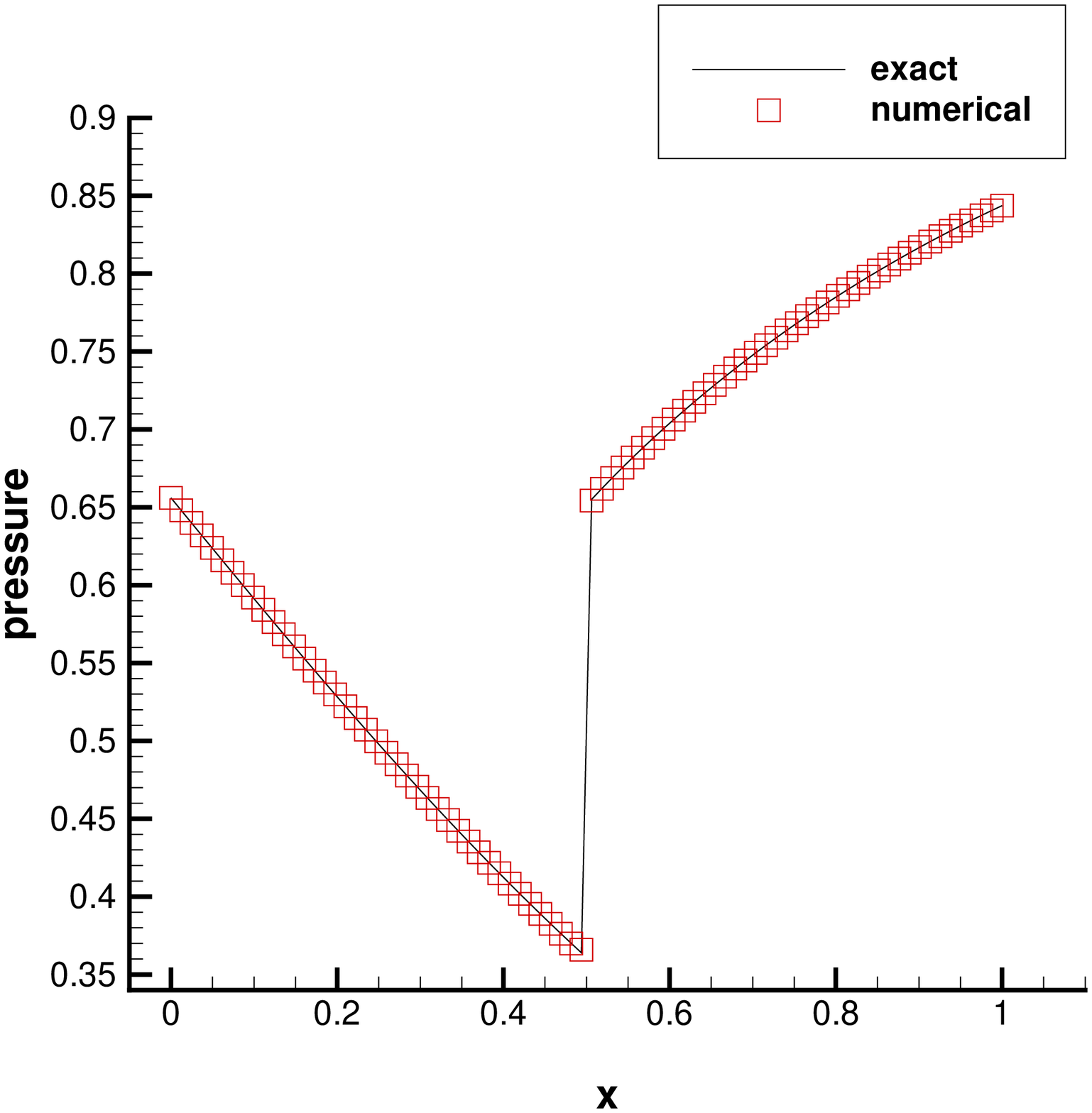}\\
\end{minipage}
\hspace{0.4cm}
\begin{minipage}[b]{0.4\textwidth}
\centering
\includegraphics[width=6.5cm]{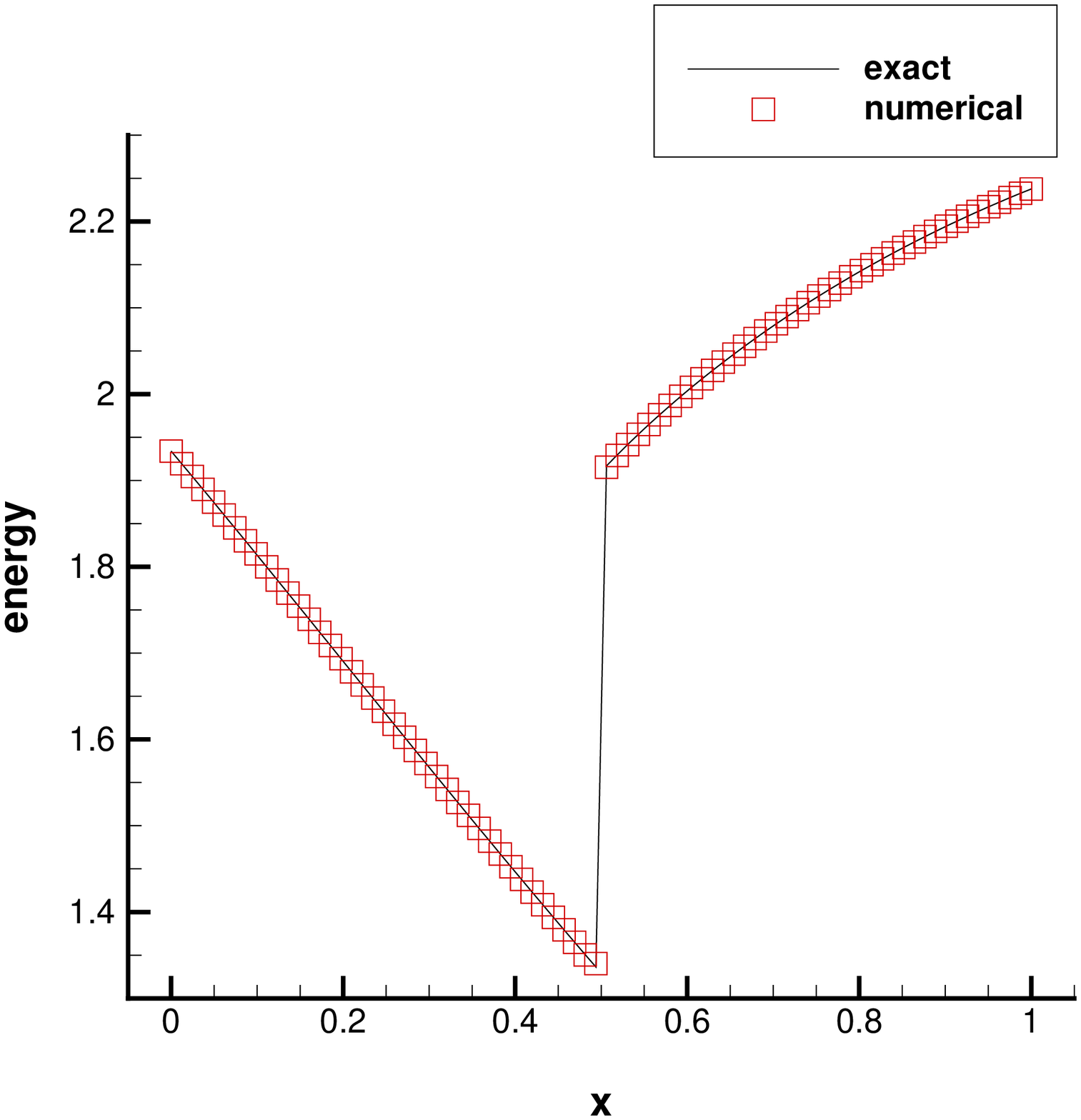}\\
\end{minipage}
\caption{The Nozzle flow problem on uniform meshes with 81 cells. Solid lines: exact solution; symbols: numerical solution. Top left: density; top right: momentum; bottom left: pressure; bottom right: total energy.}\label{Fig:nozzle1d}
\end{figure}

\begin{figure}[ht!]
\centering
\begin{minipage}[b]{0.4\textwidth}
\centering
\includegraphics[width=6.5cm]{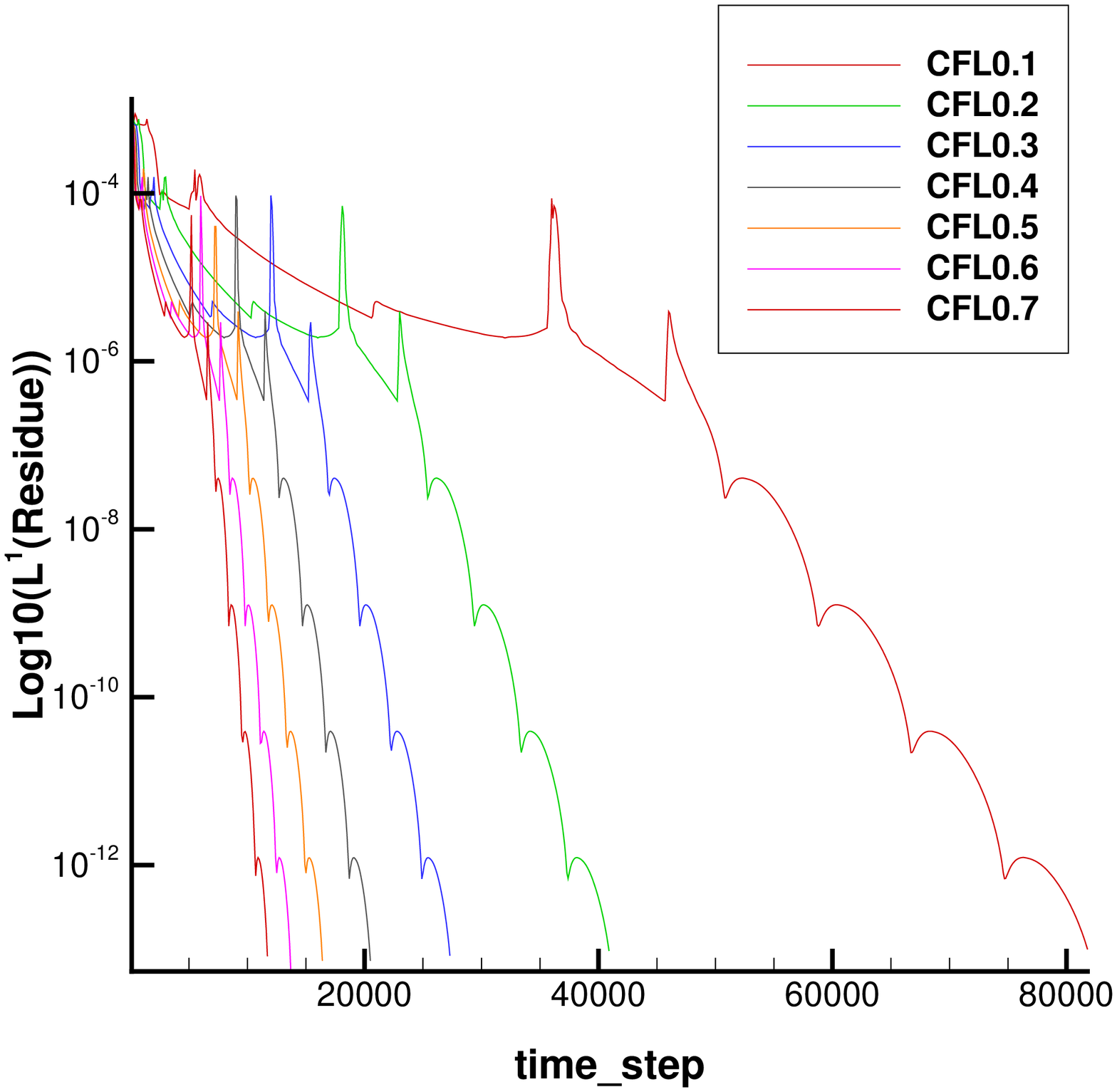}\\
\end{minipage}
\hspace{0.4cm}
\begin{minipage}[b]{0.4\textwidth}
\centering
\includegraphics[width=6.5cm]{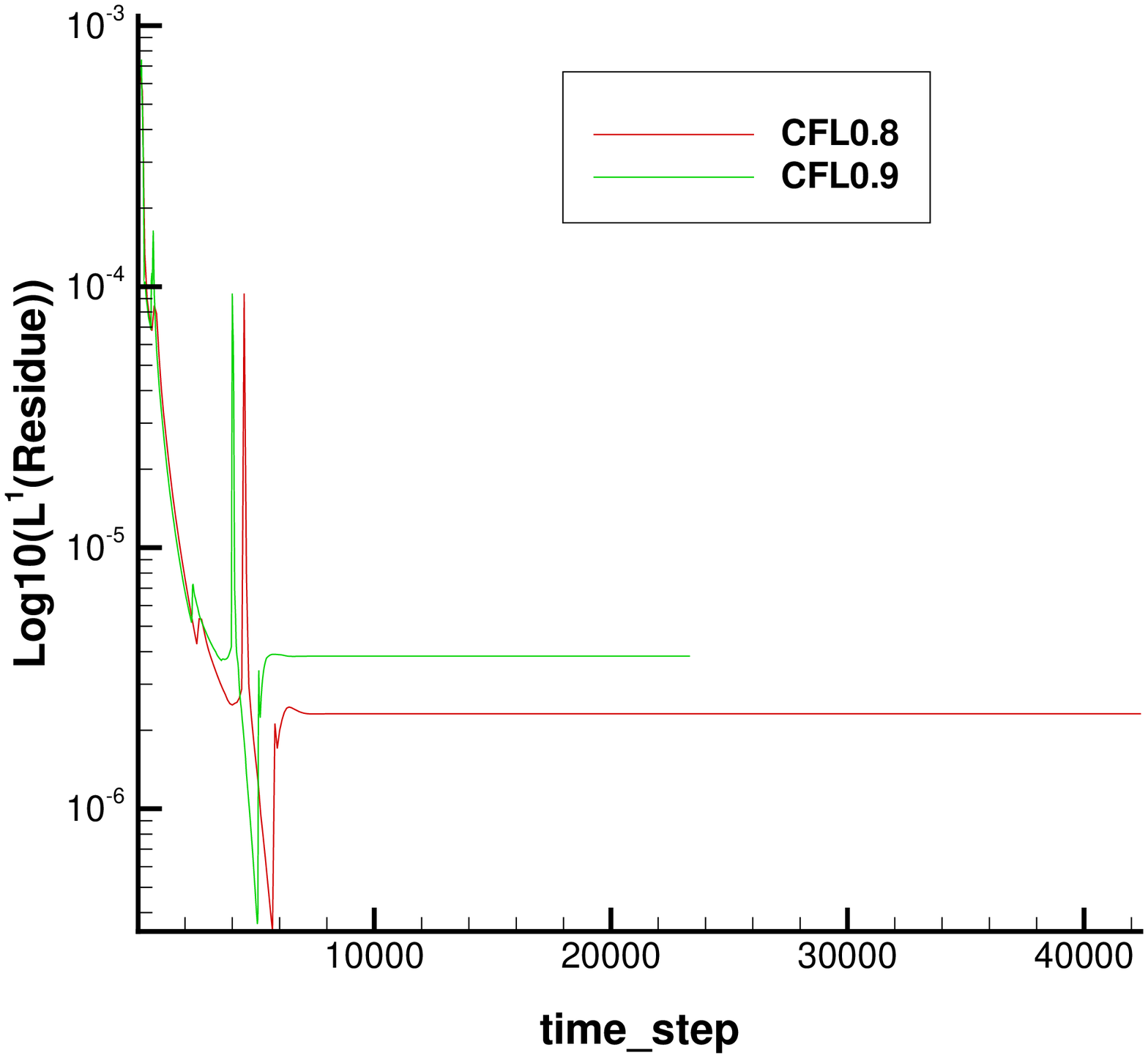}\\
\end{minipage}

\caption{The convergence histories of $L^{1}$ residue for the Nozzle flow problem}\label{Fig:nozzle_conv}
\end{figure}

\end{exmp}

\subsection{The two-dimensional scalar problems}

In this section, the numerical steady state is obtained with $L^{1}$ residue reduced to the round-off level.
\begin{exmp}\label{Eg:bur2d_sm}
We solve the steady state problem of two-dimensional Burgers equation with a source term
\begin{equation}
u_{t} + \left(\frac{1}{\sqrt{2}}\frac{u^{2}}{2}\right)_{x} + \left(\frac{1}{\sqrt{2}}\frac{u^{2}}{2}\right)_{y} = \sin\left(\frac{x+y}{\sqrt{2}}\right) \cos\left(\frac{x+y}{\sqrt{2}}\right),
\end{equation}
where $(x, y) \in \left[0, \frac{\pi}{\sqrt{2}}\right]\times \left[0, \frac{\pi}{\sqrt{2}}\right]$ with the initial condition given by
\begin{equation}
u(x, y, 0) = \beta \sin\left(\frac{x+y}{\sqrt{2}}\right).
\end{equation}
This is the one-dimensional problem studied in the \Cref{Eg:bur1d_sm} along the northeast-southwest diagonal. Since our grids are not aligned with the diagonal, this is a truly two-dimensional test case. Here we take the boundary conditions to be the exact solution of the steady state problem.

For this example, we take $\beta = 1.2$, which gives a smooth steady state solution $u(x, y, \infty) = \sin\left(\frac{x+y}{\sqrt{2}}\right)$. The errors and numerical orders are shown in the \Cref{Tab:bur2d_sm}. It can be seen clearly that the fourth order accuracy is achieved.

\begin{table}[ht!]\small
\caption{Errors and numerical orders of accuracy for the fourth order SUPG-like RD finite difference WENO-ZQ scheme for the \Cref{Eg:bur2d_sm} on uniform meshes with $N\times N$ cells}
\vspace{0.1 in}
\centering
\begin{tabular}{lllll}
\hline
\text{$N\times N$} & \text{$L^{1}$ error} & \text{Order} & \text{$L^{\infty}$ error} & \text{Order} \\
\hline
     $20\times 20$ &     7.35E-06 & &     4.29E-06 & \\ 
     $40\times 40$ &     5.61E-07 &     3.71 &     2.85E-07 &     3.91 \\ 
     $80\times 80$ &     3.86E-08 &     3.86 &     1.81E-08 &     3.98 \\ 
    $160\times 160$ &     2.53E-09 &     3.93 &     1.13E-09 &     3.99 \\
    $320\times 320$ &     1.62E-10 &     3.97 &     7.09E-11 &     4.00 \\ \hline
\end{tabular}
\label{Tab:bur2d_sm}
\end{table}

\end{exmp}

\begin{exmp}\label{Eg:bur2d_nsm}
We consider the {\color{red}steady state} solution of the following problem:
\begin{equation}
u_{t} + \left(\frac{1}{\sqrt{2}}\frac{u^{2}}{2}\right)_{x} + \left(\frac{1}{\sqrt{2}}\frac{u^{2}}{2}\right)_{y} = -\pi \cos(\pi \frac{x+y}{\sqrt{2}})u,
\end{equation}
where $(x, y)\in \left[0, \frac{1}{\sqrt{2}}\right]\times
\left[0, \frac{1}{\sqrt{2}}\right]$. This is the one-dimensional problem in the \Cref{Eg:bur1d_nsm2} along the northeast-southwest diagonal line. Inflow boundary conditions are given by the exact solution of the steady state problem. Again, since our grids are not aligned with the diagonal line, this is a truly two-dimensional test case. As before, this problem has two steady state solutions with shocks
\[
u(x , y, \infty) =\begin{cases}
1- \sin\left(\pi \frac{x+y}{\sqrt{2}}\right)  & \text{if $0 \leq \frac{x+y}{\sqrt{2}} < x_{s}$}, \\
-0.1 - \sin\left(\pi \frac{x+y}{\sqrt{2}}\right) & \text{if $x_{s} \leq \frac{x+y}{\sqrt{2}} \leq 1$},
\end{cases}
\]
where $x_s = 0.1486$ or $x_{s} = 0.8514$. Both solutions satisfy the Rankine-Hugoniot jump condition and the entropy conditions, but only the one with the shock at $\frac{x+y}{\sqrt{2}}= 0.1486$ is stable for a small perturbation.

The initial condition is given by
\[
u(x , y, 0) = \begin{cases}
\hfil 1 & \text{if $0\leq \frac{x+y}{\sqrt{2}} < 0.5$}, \\
\hfil -0.1 & \text{if $0.5 \leq \frac{x+y}{\sqrt{2}} \leq 1$},
\end{cases}
\]
where the initial jump is located in the middle of the positions of the shocks in the two admissible steady state solutions. From the \Cref{Fig:bur2d_nsm}, we can see the correct shock location and a good resolution of the solution. {\color{red} We also observe the convergence histories by different CFL numbers and the results are shown in \Cref{Fig:bur2d_nsm_conv}. We can see that the CFL number influences the convergence history, the larger CFL number and the faster convergence. }

\begin{figure}[ht!]
\centering
\begin{minipage}[b]{0.4\textwidth}
\centering
\includegraphics[width=6.5cm]{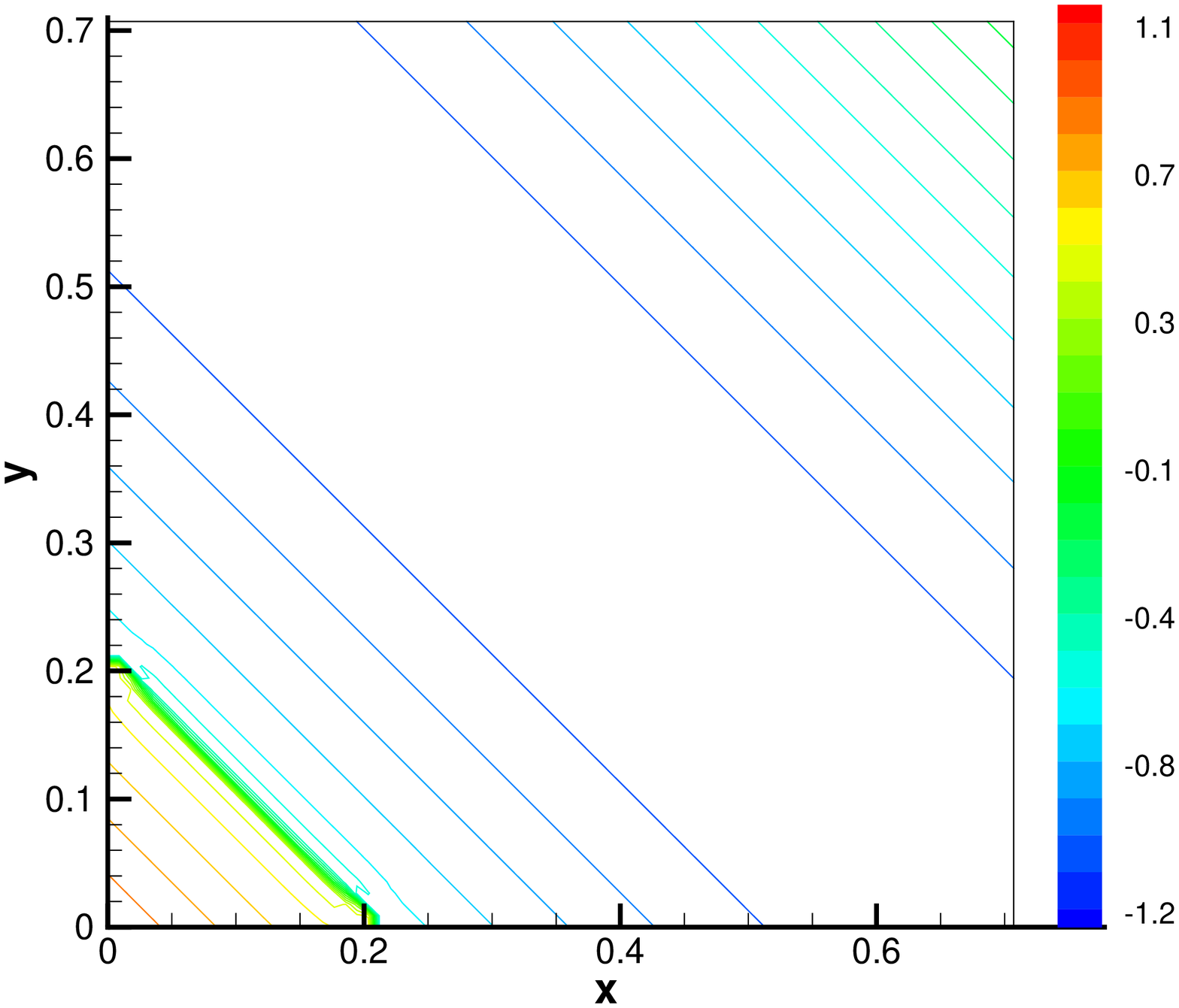}\\
\end{minipage}
\hspace{0.4cm}
\begin{minipage}[b]{0.4\textwidth}
\centering
\includegraphics[width=6.5cm]{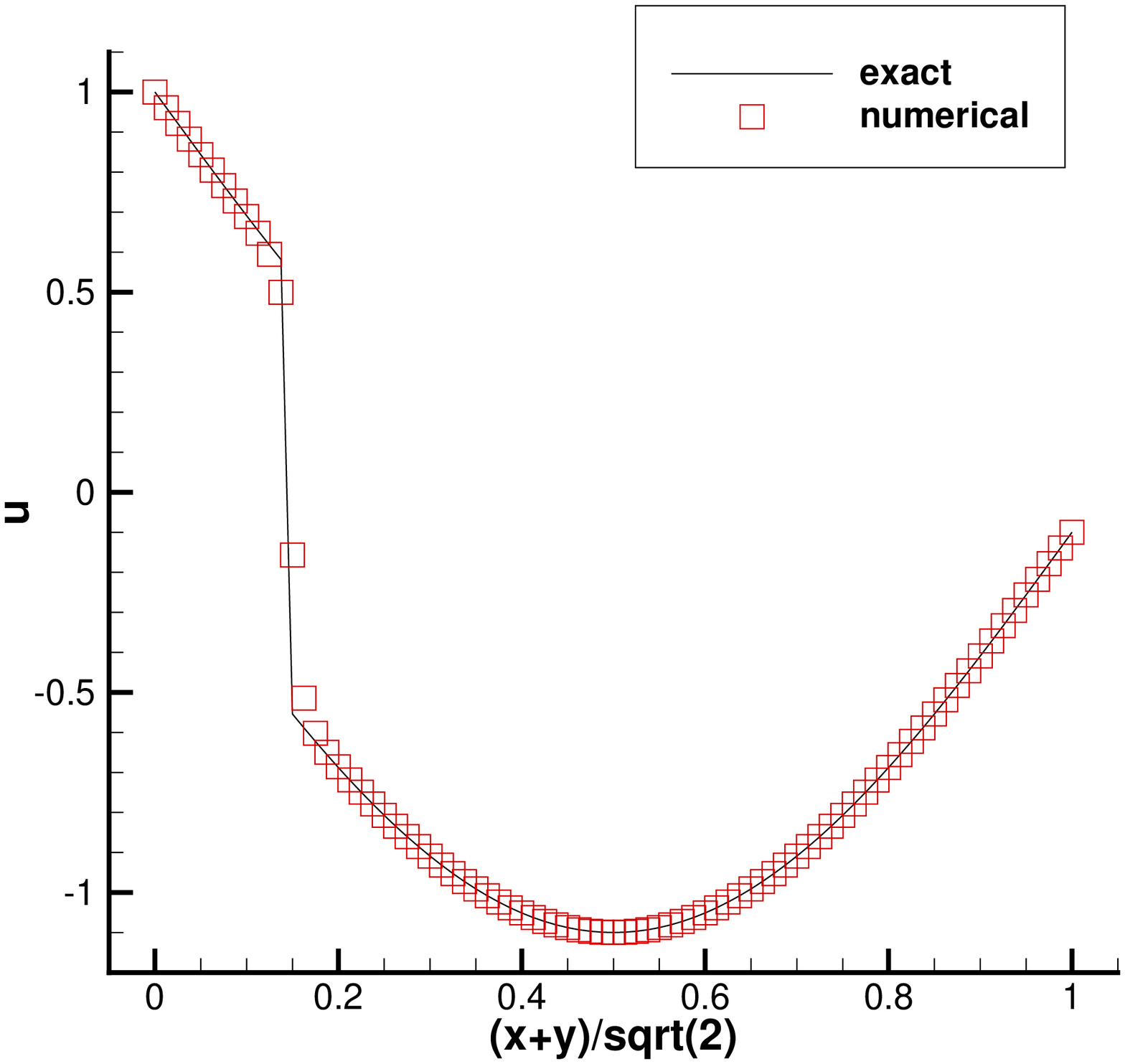}\\
\end{minipage}

\caption{The \Cref{Eg:bur2d_nsm} on uniform meshes with $80\times 80$ cells. Left: 25 equally spaced contours of the solution from -1.2 to 1.1; right: the numerical solution (symbols) versus the exact solution (solid line) along the cross-section through the northeast to southwest diagonal.}\label{Fig:bur2d_nsm}
\end{figure}

\begin{figure}[ht!]
\centering
\begin{minipage}[b]{0.4\textwidth}
\centering
\includegraphics[width=6.5cm]{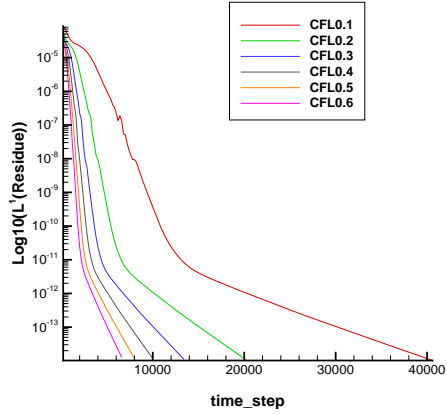}\\
\end{minipage}

\caption{The convergence histories of $L^{1}$ residue for the \Cref{Eg:bur2d_nsm}}\label{Fig:bur2d_nsm_conv}
\end{figure}
\end{exmp}

\begin{exmp}\label{Eg:bur2d_nsm2}
We consider the one-dimensional Burgers equation viewed as a two-dimensional steady state problem
\begin{equation}
u_{t} + \left(\frac{u^{2}}{2}\right)_{x} + u_{y} = 0, ~~(x, y)\in[0, 1]\times[0, 1]
\end{equation}
with the boundary conditions
\begin{equation*}
u(x, 0, t) = 1.5-2x, ~~u(0, y, t) = 1.5, ~~u(1, y, t) = -0.5.
\end{equation*}
The exact solution consists in a fan that merges into a shock which foot is located at $(x, y) = \left(\frac{3}{4}, \frac{1}{2}\right)$. More precisely, the exact solution is 
\[
u(x, y) = \begin{cases}
\text{if $y\geq 0.5$} & \begin{cases}
-0.5 & \text{if $-2(x - 3/4) + (y -1/2) \leq 0$}, \\
1.5 & \text{else}, 
\end{cases} \\
\hfil \text{else} & \text{$\max \left(-0.5, \min  \left(1.5, \frac{x-3/4}{y-1/2}\right)\right)$}.
\end{cases}
\]

This problem was studied in \citep{Cai.Gottlieb.Shu_MC1989} as a prototype example for shock boundary layer interaction. The initial condition is taken to be $u(x, y, 0) = u(x, 0, 0) = 1.5 - 2x$. The isolines of the numerical solution and the cross-sections for $y = 0.25$ across the fan, for $y = 0. 5$ right at the junction where the fan becomes a single shock, and at $y = 0.75$ across the shock, are displayed in the \Cref{Fig:bur2d_nsm2}. We can clearly observe good resolution of the numerical scheme for this example.

\begin{figure}[ht!]
\centering
\begin{minipage}[b]{0.4\textwidth}
\centering
\includegraphics[width=6.5cm]{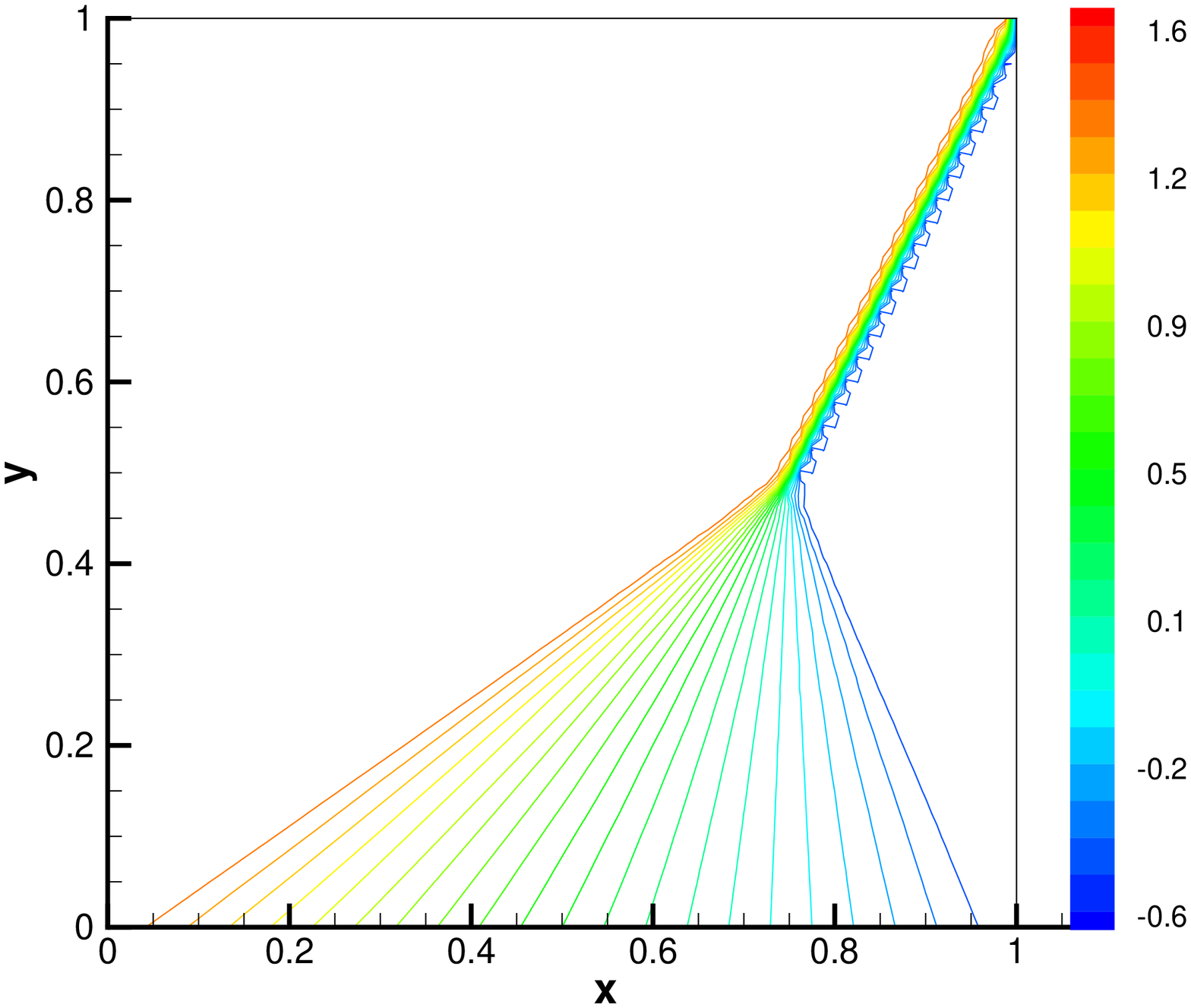}\\
\end{minipage}
\hspace{0.4 cm} 
\begin{minipage}[b]{0.4\textwidth}
\centering
\includegraphics[width=6.5cm]{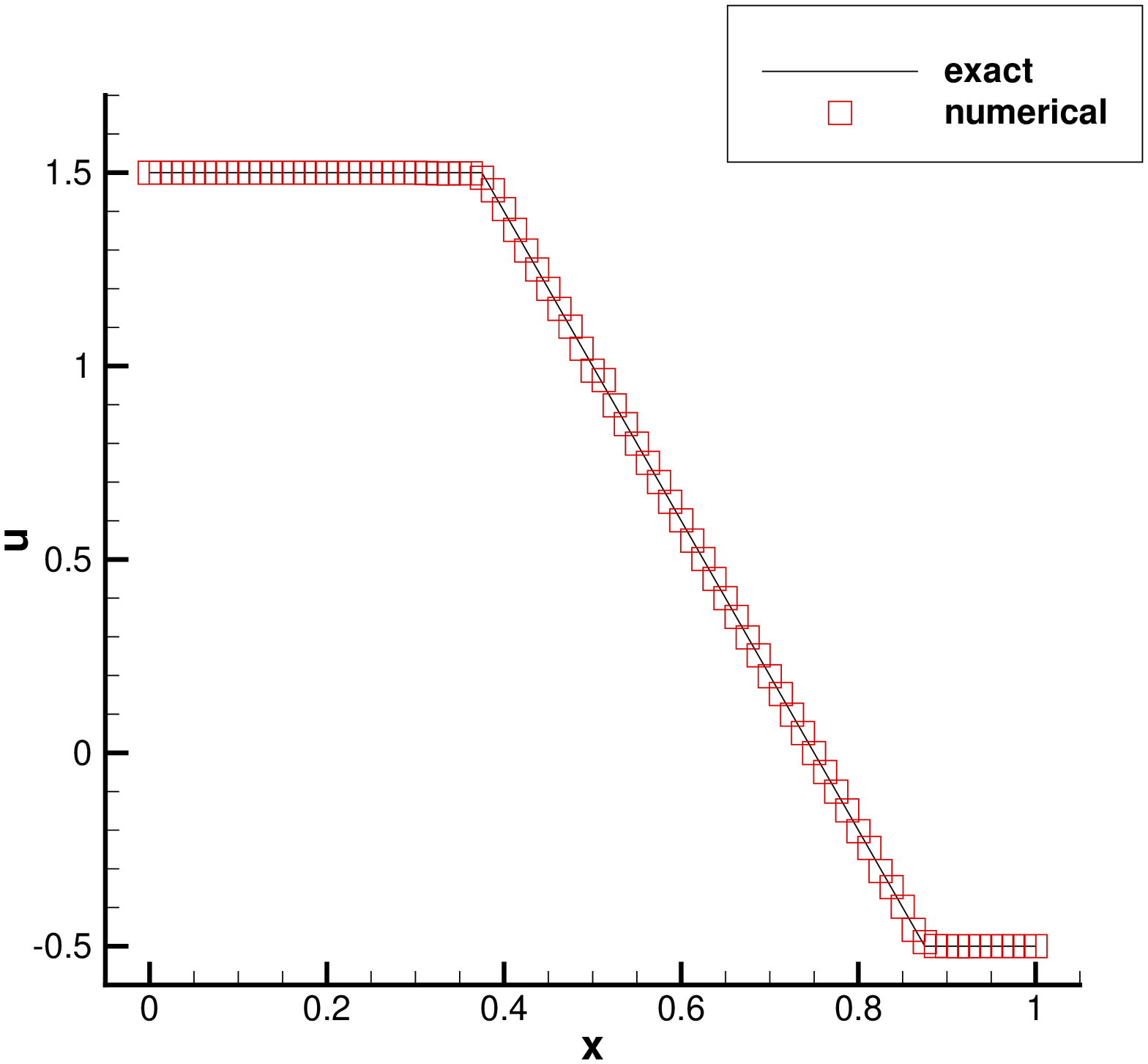}\\
\end{minipage}

\begin{minipage}[b]{0.4\textwidth}
\centering
\includegraphics[width=6.5cm]{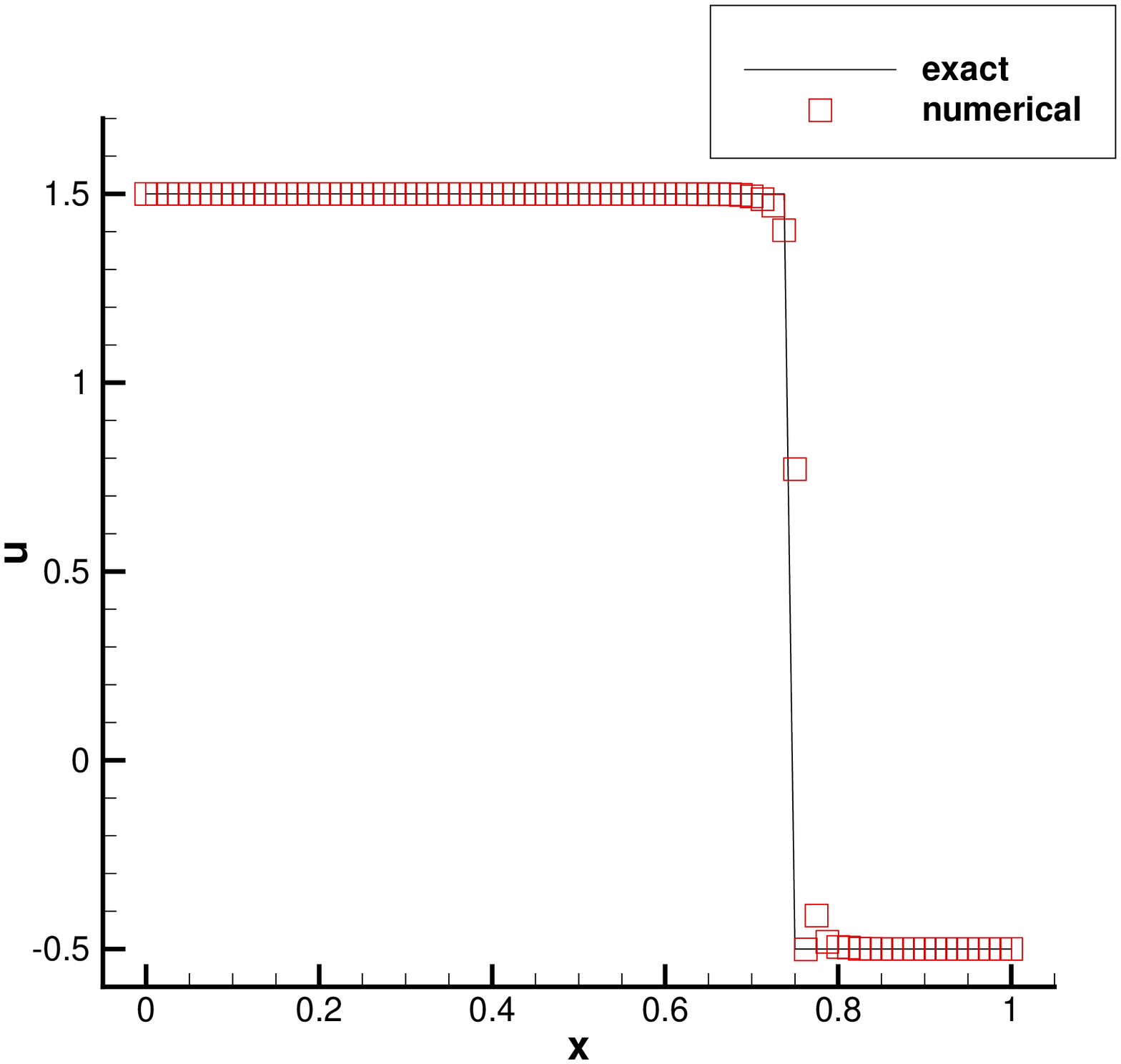}\\
\end{minipage}
\hspace{0.4 cm}
\begin{minipage}[b]{0.4\textwidth}
\centering
\includegraphics[width=6.5cm]{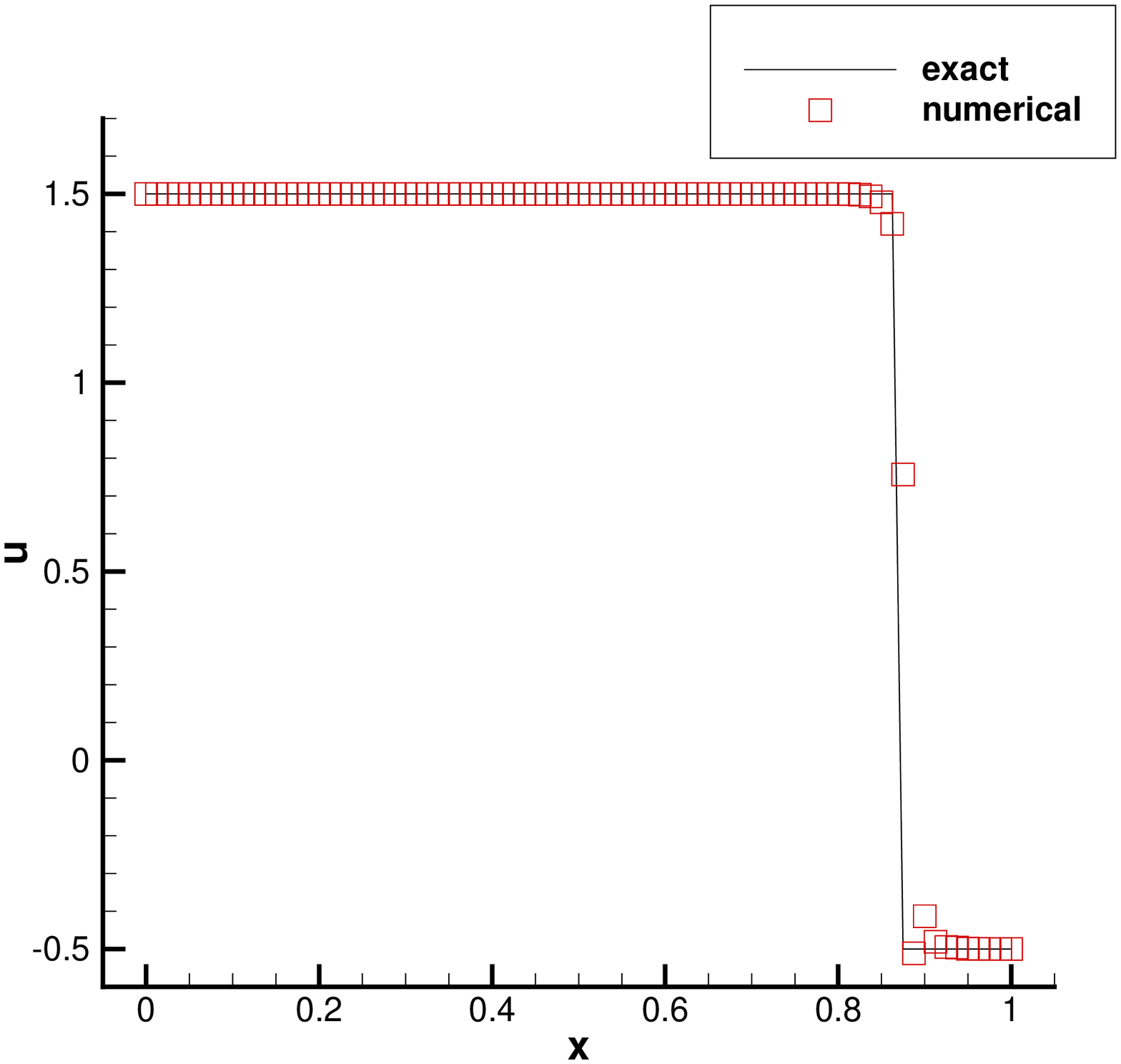}\\
\end{minipage}
\caption{The \Cref{Eg:bur2d_nsm2} on uniform meshes with $80 \times 80$ cells. Top left: 25 equally spaced contour lines from -0.6 to 1.6. Top right: cross section at $y = 0.25$; bottom left: cross section $y = 0.5$; bottom right: cross section at $y = 0.75$. For the cross section, the solid lines are for the exact solution and symbols are for the numerical solution.}\label{Fig:bur2d_nsm2}
\end{figure}

%
\end{exmp}

\subsection{The two-dimensional systems}
\begin{exmp}\label{Eg:CRP}
We consider a Cauchy-Riemann problem
\begin{equation}
\frac{\partial W}{\partial t} + A\frac{\partial W}{\partial x} + B\frac{\partial W}{\partial y} = 0, ~~ (x, y)\in [-2, 2]\times[-2, 2], ~~ t > 0,
\end{equation}
where
\begin{equation}\label{CRP_Jac}
A = \left(\begin{array}{cc}
1 & 0 \\
0 & -1
\end{array}\right)~~\text{and}~~ B = \left(\begin{array}{cc}
0 & 1 \\
1 & 0
\end{array} \right)
\end{equation}
with the following Riemann data $W = (u, v)^{T}$:
\begin{equation}\label{CRP_RMdata}
u = \begin{cases}
\hfil 1 & \text{if $x>0$ and $y>0$} \\
\hfil -1 & \text{if $x<0$ and $y>0$} \\
\hfil -1 & \text{if $x<0$ and $y<0$} \\
\hfil 1 & \text{if $x<0$ and $ y<0$}
\end{cases} ~~\text{and}~~ v = \begin{cases}
\hfil 1 & \text{if $x>0$ and $y>0$} \\
\hfil -1 & \text{if $x<0$ and $y>0$} \\
\hfil -1 & \text{if $x>0$ and $y<0$} \\
\hfil 2 & \text{if $x<0$ and $y<0$}
\end{cases}.
\end{equation}

The solution is self-similar, and therefore $W(x, y, t) = \tilde{W}\left(\frac{x}{t},\frac{y}{t}\right)$. Let $\xi = \frac{x}{t}$, $\eta = \frac{y}{t}$, then $\tilde{W}$ satisfies
\begin{equation}
(-\xi I + A)\frac{\partial\tilde{W}}{\partial\xi} + (-\eta I + B)\frac{\partial \tilde{W}}{\partial \eta} = 0,
\end{equation}
which can be written as 
\begin{equation}\label{Eq:CRPs}
\frac{\partial}{\partial \xi}[(-\xi I + A)\tilde{W}] + \frac{\partial}{\partial \eta}[(-\eta I + B)\tilde{W}] =-2\tilde{W}
\end{equation}
with the boundary conditions at infinity given by the Riemann data in \eqref{CRP_Jac} and \eqref{CRP_RMdata} at time $t = 1$. \Cref{Eq:CRPs} can be solved by RD method with boundary conditions set as the exact solution and the same initial condition as in \eqref{CRP_exact}.
\begin{align}\label{CRP_exact}
u = \begin{cases}
\hfil 1 & \text{if $x>1$ and $y>1$} \\
\hfil -1 & \text{if $x>1$ and $y<1$} \\
\hfil -1 & \text{if $x<1$ and $y>1$} \\
\hfil 1.5 & \text{if $x<1$ and $-1<y<1$}\\
\hfil 1 & \text{if $x<1$ and $y<-1$}
\end{cases} ~~\text{and}~~ v = \begin{cases}\hfil
\hfil 1 & \text{if $x>-1$ and $y>1$}\\
\hfil -1 & \text{if $x<-1$ and $y<1$} \\
\hfil -1 & \text{if $x>-1$ and $y<1$} \\
\hfil 1.5 & \text{if $x<-1$ and $-1< y <1$} \\
\hfil 2 & \text{if $x<-1$ and $y<-1$}
\end{cases}.
\end{align}
The numerical results are shown in the \Cref{Fig:crp2d}. {\color{red} From \Cref{Fig:crp_conv}, we can see $L^{1}$ residue stagnates at $10^{-6}$ level.}
\begin{figure}[ht!]
\centering
\begin{minipage}[b]{0.4\textwidth}
\centering
\includegraphics[width=6.5 cm]{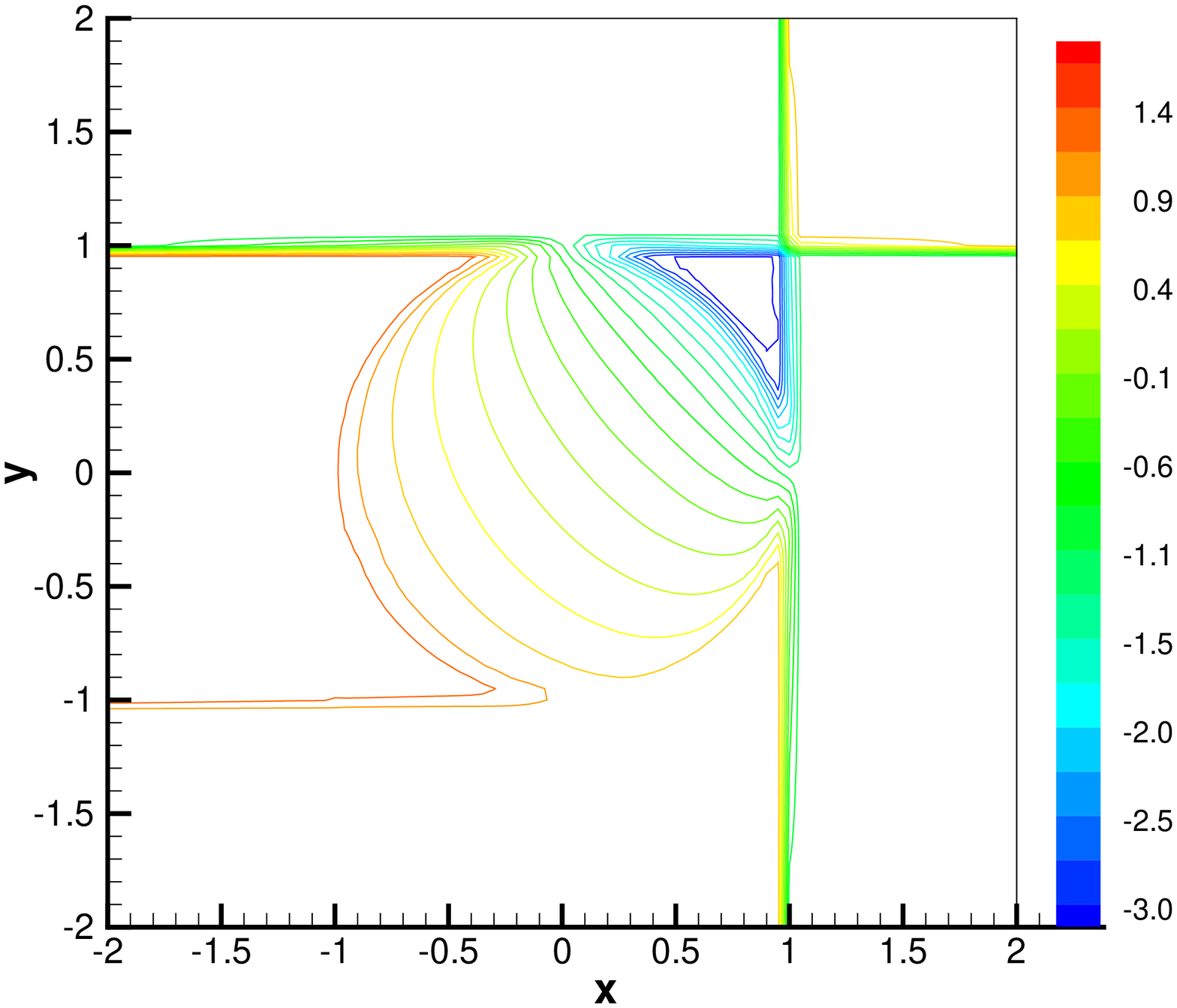}\\
\end{minipage}
\hspace{0.4 cm}
\begin{minipage}[b]{0.4\textwidth}
\centering
\includegraphics[width=6.5 cm]{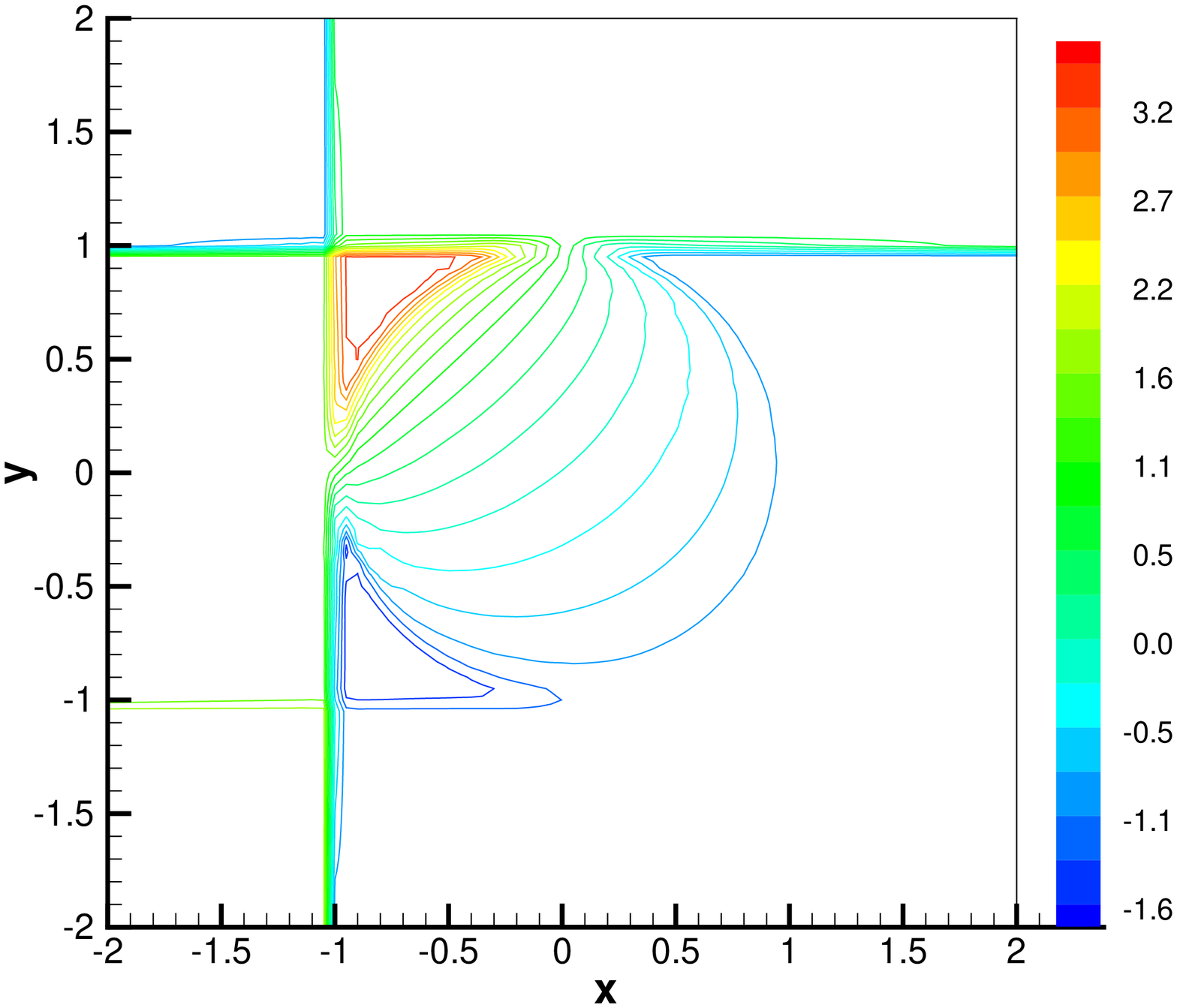}\\
\end{minipage}

\caption{The \Cref{Eg:CRP} on uniform meshes with $80\times 80$ cell. 20 Equally spaced contours for u from -3 to 1.6 (left) and 20 equally spaced contour for v from -1.6 to 3.5 (right)}\label{Fig:crp2d}
\end{figure}

\begin{figure}[ht!]
\centering
\begin{minipage}[b]{0.4\textwidth}
\centering
\includegraphics[width=6.5cm]{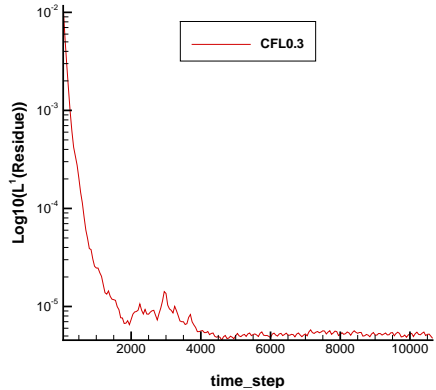}\\
\end{minipage}

\caption{The convergence history of $L^{1}$ residue for the \Cref{Eg:CRP}}\label{Fig:crp_conv}
\end{figure}
\end{exmp}

\begin{exmp}
We consider a regular shock reflection problem of the steady state solution of the two-dimensional Euler equations
\begin{equation}
{\bf u}_{t} + {\bf f}({\bf u})_{x} + {\bf g}({\bf u})_{y} = 0, ~~(x, y)\in [0, 4]\times[0, 1],
\end{equation}
where ${\bf u} = (\rho, \rho u, \rho v, E)^{T}$, ${\bf f}({\bf u}) = (\rho u, \rho u^{2} + p, \rho uv, u(E + p))^{T}$, and ${\bf g}({\bf u}) = (\rho v, \rho uv, \rho v^{2} + p, v(E + p))^{T}$. Here $\rho$ is the density, $(u, v)$ is the velocity, $E$ is the total energy and $p = (\gamma -1)(E - \frac{1}{2}(\rho u^{2} + \rho v^{2}))$ is the pressure. $\gamma$ is the gas constant which is again taken as 1.4 in our numerical tests.

The initial condition is taken to be 
\[
(\rho, u, v, p) = \begin{cases}
(1.69997, 2.61934, -0.50632, 1.52819) & \text{on $y = 1$},\\
\hfil (1, 2.9, 0, \frac{1}{\gamma}) & \text{otherwise}.
\end{cases}
\]
The boundary conditions are given by
\[
(\rho, u, v, p) = (1.69997, 2.61934, -0.50632, 1.52819) ~~ \text{on $y = 1$}
\]
and reflective boundary condition on $y = 0$. The left boundary at $x = 0$ is set as inflow with $(\rho, u, v, p)=(1, 2.9, 0, \frac{1}{\gamma})$, and the right boundary at $x = 4$ is set to be an outflow with no boundary conditions prescribed. The numerical results are shown in the \Cref{Fig:shockref}. We can clearly see a good resolution of the incident and reflected shocks. {\color{red} From \Cref{Fig:shockref_conv}, we can see the $L^{1}$ residue stagnates at $10^{-4}$  level.}
\begin{figure}[ht!]
\centering
\begin{minipage}[b]{0.4\textwidth}
\centering
\includegraphics[width=6.5cm]{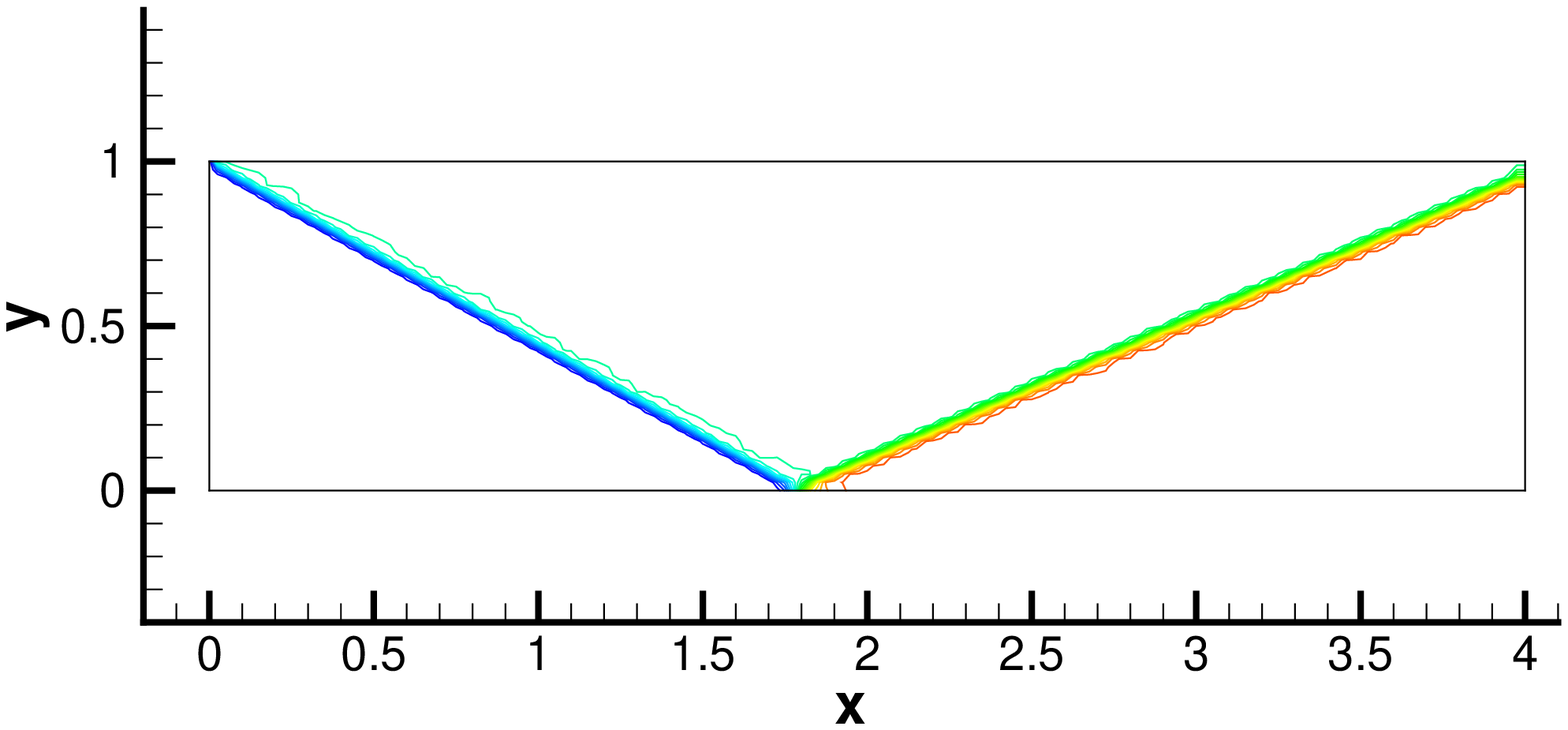}\\
\end{minipage}
\hspace{0.4cm}
\begin{minipage}[b]{0.4\textwidth}
\centering
\includegraphics[width=6.5cm]{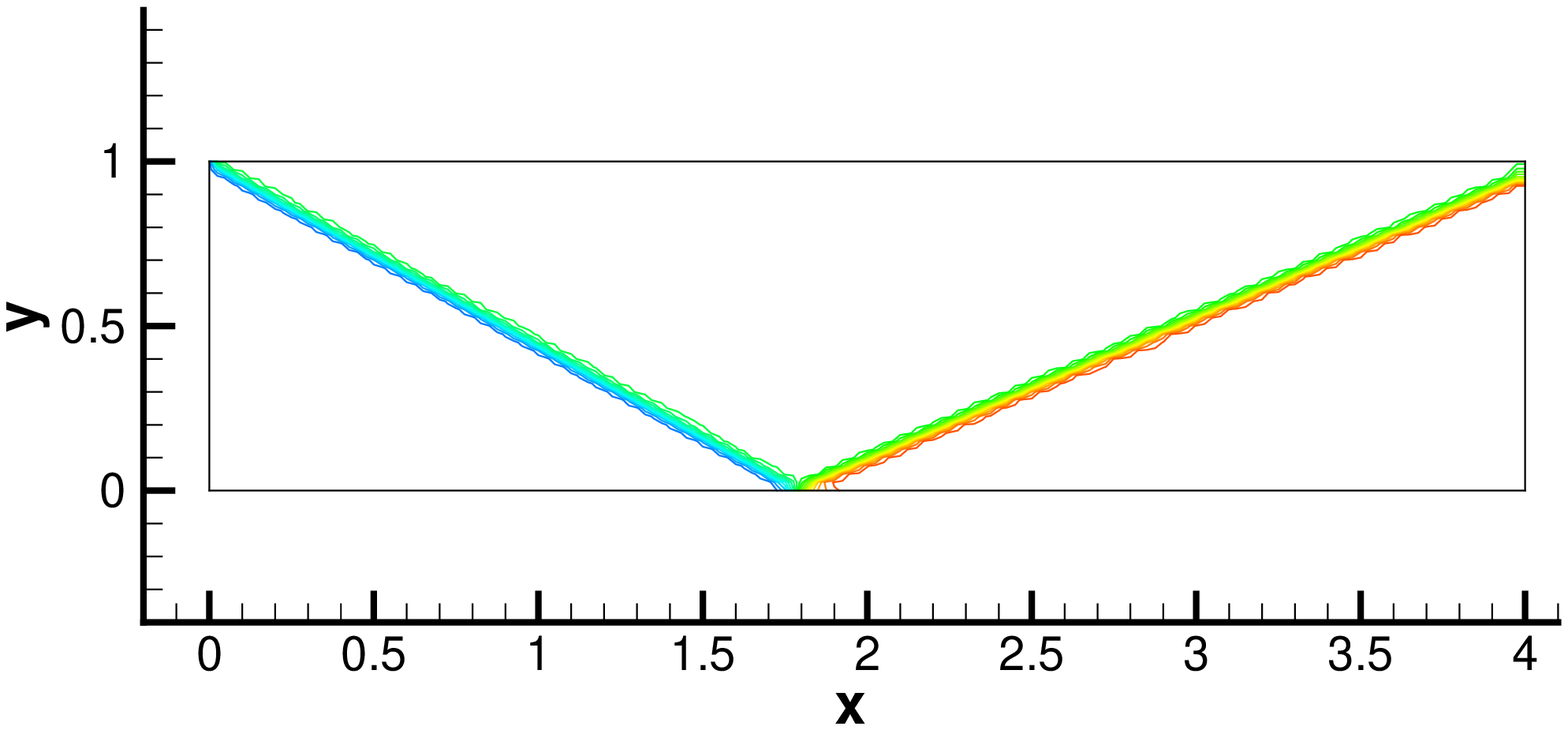}\\
\end{minipage}

\caption{Shock reflection on $160\times 40$ uniform meshes. Left: 23 equally spaced contours from 0.94 to 2.72 for the density; right: 25 equally spaced contours from 5 to 15.2 for the energy.}\label{Fig:shockref}
\end{figure}

\begin{figure}[ht!]
\centering
\begin{minipage}[b]{0.4\textwidth}
\centering
\includegraphics[width=6.5cm]{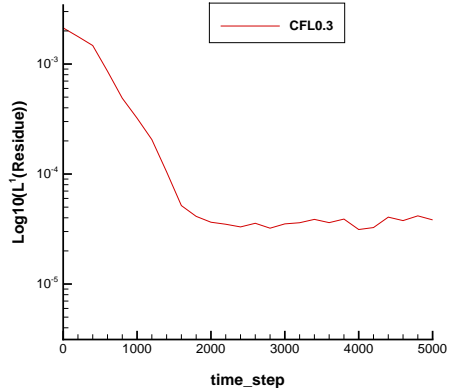}\\
\end{minipage}

\caption{The convergence history of $L^{1}$ residue for the shock reflection}\label{Fig:shockref_conv}
\end{figure}
\end{exmp}

\section{Concluding remarks}

In this paper, we proposed high order residual distribution conservative finite difference WENO-ZQ scheme for solving steady state hyperbolic equations with source terms on uniform meshes. The method is based on the WENO-ZQ integration reconstruction to achieve high order accuracy. The idea of residual distribution is adapted and allows us to obtain high order accuracy for steady state problems. We applied this proposed method to both scalar and system test problems including Burgers equation, shallow water equations, nozzle flow problem, Cauchy Riemann problem and Euler equations. In all simulations, we observed that we get the fourth order in smooth cases, and clearly see the high resolution around a shock. Future work includes using triangle meshes and extend to unsteady problems.
\section*{Acknowledgments}

J. Lin and J. Qiu are partly supported by  NSFC grant 11571290 and NSAF grant U1630247, J. Lin also is supported by the China Scholarship Council and SNF grant FZEB-0-166980. This work was performed while the first author was visiting the Institute of Mathematics, University of Zurich.


%

\section*{References}

\bibliographystyle{plainnat}
\bibliography{biblib.bib}

\end{document}